\documentstyle[11pt,amsfonts,amssymb]{article}
\def\cz{{\Bbb C}}
\def\rz{{\Bbb R}}

\def\gz{{\Bbb Z}}
\def\pr{{\rm I \! P}}
\def\qz{{\Bbb Q}}
\newcommand{\skp}{\hspace{1pt}}
\newcommand{\eps}{\varepsilon}
\newcommand{\ph}{\varphi}
\newcommand{\tr}{{\rm tr\skp}}
\newcommand{\id}{{\mbox{\rm Id}\skp}}

\newcommand{\reg}{{\mbox{\rm\scriptsize reg}}}

\newcommand{\can}{{\mbox{\rm\scriptsize can}}}

\newcommand{\rk}{{\mathop{\rm rk}}}

\newcommand{\coker}{{\mbox{\rm coker}\skp}}
\newcommand{\prj}{{\mbox{\rm pr}\skp}}
\newcommand{\img}{{\mbox{\rm im}\skp}}
\newcommand{\cl}{{\mbox{\rm cl}\skp}}

\newcommand{\Aut}{{\mbox{\rm Aut}\skp}}
\newcommand{\Hom}{{\mbox{\rm Hom}}}
\newcommand{\Ext}{{\mbox{\rm Ext}}}

\newcommand{\Spec}{{\mbox{\rm Spec}\skp}}

\newcommand{\supp}{{\mbox{\rm supp}\skp\skp}}

\newcommand{\ev}{{\mathop{\rm ev}\skp}}
\newcommand{\GW}{{\calg\calw}}

\newcommand{\AnOrb}{{\mbox{$\frak AnOrb$}}}
\newcommand{\Sets}{{\mbox{$\frak Sets$}}}

\newcommand{\longhookrightarrow}{{\begin{picture}(30,10)(-5,-2)
            \unitlength 1pt
            \put(0,2){\oval(4,4)[l]}
            \put(0,0){\vector(1,0){18}}
            \end{picture}}}
\newcommand{\llonghookrightarrow}{{\begin{picture}(25,6)(2.5,-4)
            \unitlength 1pt
            \put(2,2){\oval(4,4)[l]}
            \put(2,0){\vector(1,0){28}}
            \end{picture}}}
\newcommand{\llongrightarrow}{{\begin{picture}(25,6)(2.5,-4)
  \unitlength 1pt\put(0,0){\vector(1,0){30}}\end{picture}}}
\newcommand{\llongleftarrow}{{\begin{picture}(25,6)(2.5,-4)
  \unitlength 1pt\put(29,0){\vector(-1,0){30}}\end{picture}}}
\newcommand{\longdownarrow}{{\begin{picture}(0,0)
  \unitlength 1pt\put(0,8){\vector(0,-1){15}}\end{picture}}}
\newcommand{\longuparrow}{{\begin{picture}(0,0)
  \unitlength 1pt\put(0,-8){\vector(0,1){15}}\end{picture}}}

\newcommand{\arrowlbrt}{{\begin{picture}(0,0)
  \unitlength 1pt\put(-15,-15){\vector(2,1){30}}\end{picture}}}
\newcommand{\arrowltrb}{{\begin{picture}(0,0)
  \unitlength 1pt\put(-15,15){\vector(2,-1){30}}\end{picture}}}

\newcommand{\diagl}[1]%
  {\makebox[0cm]{${\scriptstyle#1\ }\longdownarrow
  \phantom{\scriptstyle#1\ }$}}
\newcommand{\diagr}[1]%
  {\makebox[0cm]{$\phantom{\ \scriptstyle#1}
  \longdownarrow{\ \scriptstyle#1}$}}
\newcommand{\updiagl}[1]%
  {\makebox[0cm]{${\scriptstyle#1\ }\longuparrow
  \phantom{\scriptstyle#1\ }$}}
\newcommand{\updiagr}[1]%
  {\makebox[0cm]{$\phantom{\ \scriptstyle#1}\longuparrow
  {\ \scriptstyle#1}$}}
\newcommand{\lbrtdiagr}[1]%
  {\makebox[0pt]{$\phantom{\ \ \ \ \atop\ \scriptstyle#1}
  \arrowlbrt{\ \atop\ \scriptstyle#1}$}}
\newcommand{\ltrbdiagr}[1]%
  {\makebox[0pt]{$\phantom{\scriptstyle#1\atop\ \ \ \ }
  \arrowltrb{\scriptstyle#1\atop\ \ \ \ }$}}
\newcommand{\di}{\partial}
\newcommand{\dbar}{{\bar\di}}

\newcommand{\calb}{{\cal B}}
\newcommand{\calc}{{\cal C}}

\newcommand{\cale}{{\cal E}}
\newcommand{\calf}{{\cal F}}
\newcommand{\calg}{{\cal G}}
\newcommand{\calh}{{\cal H}}
\newcommand{\cali}{{\cal I}}
\newcommand{\calj}{{\cal J}}
\newcommand{\calk}{{\cal K}}
\newcommand{\call}{{\cal L}}

\newcommand{\caln}{{\cal N}}
\newcommand{\calo}{{\cal O}}

\newcommand{\calt}{{\cal T}}

\newcommand{\calw}{{\cal W}}

\newtheorem{prop}{Proposition}[section]
\newtheorem{theorem}[prop]{Theorem}
\newtheorem{lemma}[prop]{Lemma}

\newtheorem{defi}[prop]{Definition}
\newcommand{\pf}{{\em Proof. }}
\newcommand{\qed}{{{\hfill$\diamond$}\vspace{1.5ex}}}

\title{Algebraic and symplectic Gromov-Witten invariants coincide}
\author{Bernd Siebert}
\date{April 15, 1998}
\begin{document}
\maketitle
\tableofcontents
\vspace{8ex}


\addcontentsline{toc}{section}{Introduction}
\noindent
{\Large\bf Introduction}
\vspace{1.5ex}

\noindent
Gromov-Witten invariants ``count'' (pseudo-) holomorphic
curves on algebraic or symplectic manifolds. This amounts
to intersection theory on moduli spaces of such curves.
Because in general these are non-compact, singular and not
of ``expected dimension'', a rigorous mathematical definition is far
from trivial. For a reasonably large class of manifolds including
Fano and Calabi-Yau manifolds this has first been done using symplectic
techniques by Ruan and Tian \cite{ruan}, \cite{ruantian}, \cite{ruantian2}.
The point in this approach is to restrict to sufficiently ``generic'' almost
complex structures (tamed by the symplectic form). Then the moduli
spaces are smooth of the expected dimension. This dimension is the
index of the Fredholm operator describing the moduli space locally
as map between appropriate Banach spaces. A certain positivity
condition has to be imposed on the manifolds to assure the existence
of a compactification by strata of lower dimensional manifolds.

The treatment of the general case required new techniques that would not
rely on genericity and replace the fundamental class of the manifold
by a homology class of the expected dimension, the {\em virtual
fundamental class}. This has first been achieved in the algebraic context
by Li and Tian by constructing a (bundle of) cone(s) inside a vector bundle
over a compactified moduli space; the virtual fundamental class is then
obtained by intersection with the zero section \cite{litian}. A similar
approach, based on Li and Tian's idea of using cones inside
vector bundles, but using the cotangent complex and
stack-theoretic language is due to Behrend and Fantechi \cite{behrend},
\cite{behrendfantechi}.

The construction of virtual fundamental classes in the symplectic
category has been carried out shortly later by a number of authors
\cite{fukaya}, \cite{litian2}, \cite{ruan2}, \cite{si1}. The basic approach
here is to write the moduli space as zero locus of a section of a finite
rank (orbi-) vector bundle over a finite dimensional manifold (or rather
orbifold). Locally this is not too hard. In references \cite{fukaya},
\cite{litian2} the crucial globalization is achieved by allowing for jumps
of dimension of the local models. This is still enough to construct a
homology class as intersection of the ``Euler class'' (of the bundle)
with the ``fundamental class'' (of the finite dimensional base). The other
two references use the author's construction of a finite rank orbibundle
over (a neighbourhood in an appropriate ambient space of) the moduli
space. This leads to a global description of the moduli space as zero locus
of a section of an orbibundle in an honest finite-dimensional orbifold.
\vspace{2ex}

It is natural to expect that for complex projective manifolds, to which
both approaches apply, algebraic and symplectic virtual fundamental
classes do agree (hence the Gromov-Witten invariants derived from
them). As both approaches follow rather different tracks this is, however,
far from obvious. The purpose of this paper is to confirm that the
expectation is indeed correct. We prove
\begin{theorem}
  Let $M$ be a complex projective manifold, $R\in H_2(M;\gz)$, $g,k\ge 0$.
  Let $\calc_{R,g,k}(M)$ be the moduli space of stable maps of genus
  $g$, with $k$ marked points and representing class $R$
  (cf.\ Section~\ref{section1.1}).
  
  Then the homology class associated to the algebraic
  virtual fundamental class (a Chow class on $\calc_{R,g,k}(M)$) as in
  \cite{behrend} and the symplectic virtual fundamental class
  as in \cite{si1} coincide.
\end{theorem}
For being specific we refer only to the constructions in \cite{behrend}
and \cite{si1}. The equivalence of the two algebraic constructions
and the four symplectic ones within their categories is another, fairly
straightforward albeit tedious matter, that we rather leave to a more
masochistic soul. For the symplectic case there are some comments
in \cite[\S3.4]{si3}.
\vspace{2ex}

Our proof has three central ingredients. First, we need
a holomorphic version of the construction of the ambient space
(denoted $\tilde Z$ in \cite[\S1.3,1.4]{si1}) into which $\calc_{R,g,k}(M)$
embeds as zero locus of a section $\tilde s$ of a finite rank orbibundle
$F$. This might not be possible globally, but we will gain analyticity up to
some smooth factor and this will be enough to make the comparison
work. The central point of the local construction is the fact that
spaces of holomorphic maps from a Riemann surface with non-empty
boundary (with continuous extension to the boundary) have a natural
structure of complex Banach manifold (Proposition~\ref{Stein_Hom}).
The rest of the local construction, leading to ``analytic Kuranishi models''
for $\calc_{R,g,k}(M)$ occupies Chapter~2.

The second part of the proof produces a cone $C(\tilde s)$ inside $F$,
as limit of the graphs of $t\cdot\Gamma_{\tilde s}$ as $t$ tends
to infinity. We will see (Chapter~3) that the fundamental classes
of the graphs also convergence, to a well-defined $(\dim \tilde
Z)$-dimensional homology class living on $C(\tilde s)$
(denoted $[C(\tilde s)]$, by abuse of notation). Moreover, if $\tilde s$
``splits off a trivial factor'' then both $C(\tilde s)$ and $[C(\tilde s)]$
split off a trivial factor too (a vector bundle with its fundamental class).
Of course, if $\tilde s$ is holomorphic the cone is holomorphic too.
It is in fact the cone that one obtains by applying \cite{behrendfantechi}
to an obstruction theory naturally coming from the differential sequence
associated to $\tilde s$. End second part.

Taken together this will be used in Chapter~4 to reduce the
comparison of algebraic and symplectic
virtual fundamental classes to a comparison of two morphisms of
two-term complexes to the cotangent complex of $\calc_{R,g,k}(M)$.
One is coming from the description as zero set of a section as just
mentioned, the other one is abstractly constructed in
\cite{behrendfantechi} from the universal family and universal
morphism to $M$. This comparison is the least obvious part of the proof.
It requires an explicit study of the abstractly defined
morphism (in derived categories) using \v Cech cochains. It is
quite satisfying to see how the $\dbar$-equation (describing
$\tilde s$) naturally arises by partial integration applied
to fiber integrals coming from the explicit version of relative (Serre-)
duality (see the proof of Lemma~\ref{integral_descr}).

The first Chapter serves two purposes. First, it contains an account
of the algebraic definition of virtual fundamental classes to fix
notations and to make the paper more self-contained. We follow
the elementary reformulation (avoiding Artin stacks) of
\cite{behrendfantechi} previously given by the author \cite{si2}.
Second we go over from algebraic to analytic spaces, or rather
from Deligne-Mumford stacks to complex analytic orbispaces.
The last Chapter~5 provides a GAGA-type result concerning
push-forwards and relative duality in algebraic and analytic
derived categories of sheaves. The main result is
Proposition~\ref{duality} which gives an explicit description of algebraic
relative duality for algebraic families of prestable curves in terms
of analytic \v Cech cochains and fiber integrals.
\vspace{2ex}

Our proof of the comparison theorem has been sketched in
some detail in the survey \cite{si3}, submitted in May 1997.
Shortly afterwards we learned from J.\ Li and G.\ Tian that
they were able to proof the same result, using their respective definitions
of algebraic and symplectic virtual fundamental classes
\cite{litian3}.


\section{Complex analytic GW-theory}
To compare algebraic and symplectic definitions of GW-invariants, as a first,
mostly trivial step, it is natural to translate the former into the category of
complex (analytic) spaces. This will be the purpose of this chapter.


\subsection{Analytic orbispaces versus Deligne-Mumford stacks}
\label{section1.1}
Given a smooth projective scheme $M$ over a field $k$ the natural arena
for GW-theory is the space $\calc(M)$ of stable marked curves
in $M$ (Kontsevich's ``stable maps''). So the$k$-rational points of
$\calc(M)$ are in one-to-one correspindence with isomorphism classes
of triples $(C,{\bf x},\ph)$ with $C$ a reduced, connected algebraic curve,
proper over $k$ and with at most ordinary double points, ${\bf x}=
(x_1,\ldots,x_k)$ a tuple of $k$-rational points in the regular locus
$C_\reg\subset C$, and $\ph:C\rightarrow M$ a $k$-morphism with
the property that
\[
  \Aut(C,{\bf x},\ph)\ =\ \{\psi\in\Aut(C)\,|\,\psi({\bf x})={\bf x}, \ph\circ\psi=\ph\}
\]
is finite ({\em stability}). With the obvious notion of families of stable
marked curves over (that is, parametrized by) $k$-schemes,
$\calc(M)$ (or rather the associated fibered groupoid) has been verified
to be a Deligne-Mumford (DM-) stack \cite{behrendmanin}.

In the analytic category, i.e.\ $k=\cz$ and $M$ viewed as complex projective
manifold, the DM-stack can be replaced by a notion of analytic orbispace
that we now introduce. This will be a generalization of both complex
orbifolds and complex spaces.

Local models for such spaces are tuples $(q: \hat U\rightarrow U, G,\alpha)$
with
\begin{itemize}
\item
  $G$ is a finite group, viewed as zero-dimensional reduced complex space
\item
  $\hat U$ is a (not necessarily reduced, but finite-dimensional) complex
  space
\item
  $\alpha: G\times \hat U \rightarrow \hat U$ is a (not necessarily effective)
  holomorphic group action on $\hat U$
\item
  $q$ is a quotient of $\hat U$ by $G$ in the category of complex spaces
  (or, equivalently, in the category of ringed spaces, cf.\ e.g.\ 
  \cite[\S49A]{kaupkaup}).
\end{itemize}
We will often just write $U= \hat U/G$ for such {\em local uniformizing
system} $(q,G,\alpha)$. The definition of analytic orbispaces now runs
completely analogous to the case of orbifolds as given in \cite{satake}: One
first defines the notion of {\em morphisms}, and in particular {\em open
embeddings\,} of local uniformizing systems. An {\em analytic orbispace
structure} on a Hausdorff space $X$ is then a covering by local uniformizing
systems $\{ U_i= \hat U_i/G_i\}_{i\in I}$ (i.e.\ $X=\bigcup U_i$) compatible on
overlaps $U_i \cap U_j \neq \emptyset$ by open embeddings. Finally, an {\em
analytic orbispace} is an equivalence class of analytic orbispace structures.
By abuse of notation we will also just write $X$ for the analytic orbispace.
Of course there is also a notion of morphisms of analytic orbispaces making
the set of analytic orbispaces into a category, denoted $\AnOrb$.

Similarly, we may introduce the notions of (topological or analytic)
{\em orbibundles} and of {\em coherent orbisheaves} on $X$ by requiring the
associated linear fiber space over $X$ (cf.\ e.g.\ \cite[\S1.4]{fischer}) to
be analytic orbispaces over $X$ with local uniformizers having a
well-defined linear structure on the fibers.

This all being a trivial translation of \cite[\S1.1]{si1}, which treats the case of
topological Banach orbifolds, to the category of complex spaces we merely
give these indications and refer to op.cit.\ for more details.
\vspace{2ex}

We claim that $\calc(M)$ has naturally the structure of an analytic
orbispace, in such a way that $\calc(M)$ represents the functor
\begin{eqnarray*}
  {\frak C}(M): \AnOrb &\longrightarrow& \Sets\\
  T&\longmapsto& \left\{
  \begin{array}{c}\mbox{analytic families of stable holomor-}\\
                             \mbox{phic curves in $M$ parametrized by $T$}
  \end{array} \right\}\Big/\mbox{iso}\,.
\end{eqnarray*}
Here we use the following
\begin{defi}\rm
Let $T$ and $M$ be analytic orbispaces. An {\em analytic family of stable
holomorphic curves in $M$} parametrized by $T$ is a tuple of morphisms
\[
  (q:X\rightarrow T, \underline{\bf x}: T\rightarrow X\times_T \ldots
  \times_T X, \Phi: X\rightarrow M)
\]
of complex orbispaces with
\begin{itemize}
\item
  $q$ is flat
\item
  for any $t\in T$, $(q^{-1}(t), \underline{\bf x} (t), \Phi|_{q^{-1}(t)})$ (with
  the induced analytic structure) is a stable marked holomorphic curve in $M$
  \cite[Def.3.5]{si1}.
  \qed
\end{itemize}
\end{defi}
\vspace{-2ex}

\noindent
The precise result is
\begin{prop}\label{CM_analyt_space}
  For any complex space $M$, ${\frak C}(M)$ is representable by
  an analytic orbi\-space $\calc(M)$.
\end{prop}
\pf
To define a local uniformizing system at some stable holomorphic curve
$(C,{\bf x}, \ph)$ in $M$ let $(q:\calc \rightarrow S, \underline{\bf x})$ be an
analytically semi-universal deformation of $(C, {\bf x})$,
cf.\ e.g.\ \cite[\S2.2]{si1}. It is well-known that $\Hom_S (\calc, M)$,
the space of morphisms from the fibers of $q$ to $M$, is a complex
space (representing the corresponding functor) \cite{douady},\cite{pourcin}.
This is almost the space we want, but if $\Aut^0 (C,{\bf x})$ is non-trivial
we have to the equivalence relation on $\Hom_S (\calc, M)$ generated by the
germ of the action of $\Aut^0 (C,{\bf x})$ on $q$. We refer to this process
as ``rigidification''.

To this end let ${\bf y}= (y_1,\ldots, y_l)$ be a tuple of points in $C$ such
that
\begin{itemize}
\item
  $(C,{\tilde{\bf x}})$ is a (Deligne-Mumford) stable curve, where
  ${\tilde{\bf x}}= (x_1,\ldots,x_k,y_1,\ldots y_l)$
\item
  for any $i$, $\ph$ is an immersion at $y_i$ (possible by stability)
\end{itemize}
The first property requires in particular finiteness of $\{\Psi\in \Aut(C,{\bf x})
\mid \Psi(y_i)=y_i\}$, so a minimal choice requires the insertion
of $l=\dim \Aut(C,{\bf x})$ points. By the second property we may choose
local Cartier divisors $H_1, \ldots,H_l \subset M$ and open disks
$U_i \subset C$ with $(\ph|_{U_i})^{-1} (H_i)
= y_i$ ($y_i$ with reduced structure). Now let
$(q: \calc \rightarrow S, \underline{\tilde{\bf x}})$ be an analytically
universal deformation of $(C,\tilde{\bf x})$ with $\tilde{\bf x}
=(x_1, \ldots, x_k, y_1,\ldots, y_l)$. Extend $U_i$ to open polycylinders (say)
$\tilde U_i \subset \calc$. By restriction to an open subspace $Z\subset
\Hom_S(\calc,M)$ we may assume that $(\psi|_{\tilde U_i})^{-1} (H_i)$
consists of exactly one reduced point for any $\psi \in Z$. Evaluation at the
deformation of the $\nu$-th marked point defines $k+l$ morphisms
\[
  \ev_\nu: Z\longrightarrow M\,.
\]
Set
\[
  \hat U:= (\ev_{k+1},\ldots,\ev_{k+l})^{-1} (H_1\times\ldots\times H_l)\,.
\]
We claim that the restriction of the universal curve $\calc\times_S
\Hom_S (\calc,M)\rightarrow \Hom_S (\calc,M)$ to $\hat U$
together with the universal (evaluation) morphism from the universal
curve to $M$, is a universal deformation of $(C,{\bf x},\ph)$. So let
$(X\rightarrow T,\underline{\bf x}',\Phi')$ be an analytic family of
stable holomorphic curves together with an isomorphism of the
fiber over some point $0\in T$ with $(C,{\bf x},\ph)$. Since $\ph$ is
transverse to $H_i$, local defining equations of $H_i$ pull back to
local holomorphic functions on $X$ that restrict to local holomorphic
coordinates for $C$ near $\ph^{-1}(H_i)$ on the central fiber. This
shows that possibly after replacing $T$ by a neighbourhood
of $0\in T$ the Cartier divisor ${\Phi'}^{-1}(H_i)$ is a section of
$X\rightarrow T$. We denote this section by $\underline{y}'_i$
and write $\underline{\tilde{\bf x}}'= (\underline{\tilde x}'_1,\ldots,
\underline{\tilde x}'_k, \underline{y}'_1,\ldots,
\underline{y}'_l)$. By the universal property of
$(\calc\rightarrow S, \underline{\tilde {\bf x}})$
there exists a unique morphism $T\rightarrow S$ such that
$(X\rightarrow T,\underline{\tilde{\bf x}}')$ is isomorphic to the
pull-back of the universal family over $S$. Moreover, there is a unique such
isomorphism inducing the given identification of $X_0$ with $C$.

In turn the universal property of the Hom-space produces a unique
morphism $T\rightarrow \Hom_S(\calc,M)$ such that $\Phi'$ is
the composition of the evaluation map and the product morphism
from $X$ to $\calc \times_S \Hom_S (\calc,M)$. And if $w_i=0$ is a
defining equation for $H_i$ then $(\Phi'\circ\underline{y}'_i)^* w_i
=0$ by definition of $\underline{y}'_i$. So the map from $T$ to
$\Hom_S(\calc,M)$ indeed factors over $\hat U$. Given any other map
from $T$ to $\hat U$ with these properties, we can pull-back the sections
$\underline{\tilde{\bf x}}$ to $T$ to see that this map coincides with the
one just constructed.

Since $\Aut(C,{\bf x},\ph)\subset \Aut(C,{\bf x})$ acts on
$(\calc\rightarrow S,\underline{\tilde{\bf x}})$ ($(\calc\rightarrow S,
\underline{\bf x})$ is a semi-universal deformation of $(C,{\bf x})$!)
the universal property now immediately implies an action of this group
on $\hat U$ and compatibility of this action with open embeddings of
local uniformizing systems, the existence of which being itself provided
by universality. A final remark concerns global existence
of the universal curve and universal morphism. In fact, the universal curve
is itself isomorphic to $\calc(M)$, the evaluation morphism is evaluation
at the last marked point and the morphism to $\calc(M)$ is by forgetting the
last marked point and stabilizing (i.e.\ successively contracting all
unstable components on which the map is trivial). These are morphisms
by work of Knudson \cite{knudson}.
\qed


\subsection{Analytic global normal space, local construction}
By the analytic analog of \cite{behrendfantechi} or its elementary
reformulation in \cite{si2} we need a (free) global normal space
(\cite[Def.3.1 and \S4.1]{si2}) for $\calc(M)$
relative ``the Artin stack of prestable curves'' ${\frak M} = \coprod_{g,k}
{\frak M}_{g,k}$ to construct the virtual fundamental class on $\calc (M)$.
Recall that a global normal space for $\calc(M)$ relative ${\frak M}$ would
be a morphism in the derived category of coherent orbisheaves on $\calc(M)$
\[
  \ph^\bullet: \calf^\bullet=[\calf^{-1}\rightarrow\calf^0]
  \longrightarrow \call^\bullet_{\calc(M)/{\frak M}}
\]
with
\begin{itemize}
\item
  $\calf^0$ and $\calf^{-1}$ are locally free
\item
  $\ph^\bullet$ induces an isomorphism in $H^0$ and an epimorphism in
  $H^{-1}$.
\end{itemize}
Here $\call^\bullet_{\calc(M)/{\frak M}}$ is the analytic cotangent complex
relative ${\frak M}$. Invoking existing literature (\cite{bingenerkosarew},
\cite{flenner}, \cite{illusie}) for the case of analytic orbispaces (or an
appropriate analog of DM-stacks in the analytic category) relative an
Artin stack is however questionable. Instead of justifying this rather
technical step we will give our construction on the level of local
uniformizing systems and show in the next section that the
corresponding quotients of cones in the complex orbibundle
$H$ with local uniformizer $F_1= L(\calf^{-1})$
globalize. It will be clear from the local construction of this cone that the
result will be exactly the complex space associated to the DM-stack of
cones in the stack-theoretic version of $H$ in \cite{behrendfantechi}.

We start with the ``obstruction theory'' for spaces of morphism as in
\cite[Ch.6]{behrendfantechi}, in the relevant relative formulation. The
following discussion is literally valid both algebraically or analytically.
For $(C,{\bf x}, \ph)\in \calc(M)$ let $(q:\calc \rightarrow S, \tilde{ \underline{
\bf x}})$ be a universal deformation of a rigidification $(C, \tilde{\bf x})$ of
$(C,{\bf x})$, $\tilde{\bf x}= (x_1, \ldots, x_k, y_1, \ldots, y_l)$, as in the
proof of Proposition~\ref{CM_analyt_space}. Over the open
subspace $Z\subset \Hom_S(\calc, M)$
lives the universal curve $\pi: \Gamma:=\calc \times_S Z \rightarrow Z$
with the universal morphism $\Phi: \Gamma \rightarrow M$.
\[\begin{array}{ccccc}
  \calc&\longleftarrow&\Gamma&\stackrel{\Phi}{\longrightarrow}M\\[2pt]
  \longdownarrow&&\diagr{\pi}\\[6pt]
  S&\longleftarrow&Z
\end{array}\]
From the functorial properties of the cotangent complex \cite{illusie},
\cite[I.2.17]{flenner} we obtain morphisms
\[
  L\Phi^*\call_M^\bullet \longrightarrow \call_\Gamma^\bullet
  \longrightarrow \call_{\Gamma/\calc}^\bullet\,,
\]
and, by flatness of $\pi$, an isomorphism \cite[I.2.26]{flenner},
\cite[II.2.2.3]{illusie}
\[
  L\pi^*\call_{Z/S}^\bullet \longrightarrow \call_{\Gamma/\calc}^\bullet\,.
\]
All these morphisms are to be understood in $D^-(\calo_\Gamma)$ the
derived category of the category of $\calo_\Gamma$-modules, bounded
to the right, or rather the derived category $D^-_{\rm coh} (\calo_\Gamma)$
of complexes with coherent cohomology. The resulting morphism $L \Phi^*
\call_M^\bullet \rightarrow L\pi^* \call_{Z/S}^\bullet$ is tensored (in the
left-derived sense) with $\omega_\pi$, the relative dualizing sheaf of $\pi:
\Gamma \rightarrow Z$. Applying $R\pi_*$ we get
\[
   R\pi_*(L\Phi^*\call_M^\bullet\stackrel{L}{\otimes}\omega_\pi)
  \longrightarrow
  R\pi_*(L\pi^*\call_{Z/S}^\bullet\stackrel{L}{\otimes}\omega_\pi)\ \simeq\ 
  \call_{Z/S}^\bullet\stackrel{L}{\otimes}R\pi_*\omega_\pi\,.
\]
Now by smoothness of $M$, $\call_M^\bullet =\Omega_M$ (viewed as complex
concentrated in degree 0), hence $L\Phi^*\call_M^\bullet\stackrel{L}{\otimes}
\omega_\pi =\Phi^* \Omega_M\otimes\omega_\pi$, and by relative duality
\cite{hartshorne}, \cite{verdier}
\[
  R\pi_*(\Phi^*\Omega_M\otimes \omega_\pi)\ \simeq\ 
  [R\pi_*\Phi^*\calt_M]^\vee,\quad
  R\pi_*\omega_\pi\ \simeq\ \calo_Z\,,
\]
where $\calt_M=\Omega_M^\vee$ is the locally free sheaf of holomorphic
vector fields on $M$. For any complex $\cale^\bullet$
of coherent sheaves, bounded to the
left, we use the notation $(\cale^\bullet)^\vee =\Ext (\cale^\bullet,
\calo)$ for the dual complex in the derived sense.
In particular, $\calf^\vee= \Hom_\calo (\calf,
\calo)$ for locally free coherent sheaves. We have thus produced
a morphism (which is in fact an ``obstruction theory for $Z$ relative $S$''
in the notation of \cite{behrendfantechi})
\begin{eqnarray}
  [R\pi_*\Phi^* \calt_M]^\vee &\longrightarrow& \call_{Z/S}^\bullet\,.
  \label{obstr_theory_I}
\end{eqnarray}
To represent the left-hand side by a morphism of locally free sheaves one
needs to assume $M$ projective. Let $\calh_M$ be an ample line bundle on
$M$. Then by stability
\[
  \call\ :=\ \omega_\pi (\underline x_1\ldots+\underline x_k) \otimes
  \Phi^*\calh_M^{\otimes 3}
\]
has positive degree on each irreducible component of any fiber of $\pi$, where
we wrote $\underline x_i$ for the Cartier divisor corresponding
to the $i$-th entry of $\underline{\tilde {\bf x}}$. It is not hard to see that
choosing $\nu$ large enough the natural morphism
\[
  \caln\ :=\ \pi^*\pi_*(\Phi^*\calt_M \otimes\call^{\otimes \nu}) \otimes
  \call^{\otimes -\nu} \longrightarrow \Phi^*\calt_M
\]
is surjective and $\pi_*\caln =\pi_*\calk =0$, $\calk:= \ker (\caln\rightarrow
\Phi^*\calt_M)$ (cf.\ \cite[Prop.5]{behrend}). Pushing forward the exact sequence
\begin{eqnarray}\label{twist_sequence}
  0\longrightarrow \calk \longrightarrow \caln \longrightarrow \Phi^*\calt_M
  \longrightarrow 0
\end{eqnarray}
by $\pi$ in the derived sense we obtain an exact triangle
\[
  R\pi_* \Phi^*\calt_M[-1] \longrightarrow [0\rightarrow R^1\pi_*\calk]
  \longrightarrow [0\rightarrow R^1\pi_*\caln]\longrightarrow
  R\pi_*\Phi^*\calt_M\,,
\]
hence an isomorphism of the mapping cone of the middle morphism with
$R\pi_* \Phi^*\calt_M$:
\[
  [R^1\pi_*\calk\rightarrow R^1\pi_*\caln] \stackrel{\simeq}{\longrightarrow}
  R\pi_* \Phi^*\calt_M\,.
\]
Note that $R^1\pi_*\calk$, $R^1\pi_*\caln$ are locally free for $\pi_*\calk=
\pi_*\caln =0$, The derived dual of $[R^1\pi_*\calk \rightarrow R^1\pi_*
\caln]$ may thus be taken entrywise. Write $\calg := (R^1\pi_*\calk)^\vee$,
$\calh:= (R^1\pi_*\caln)^\vee$, and $G$, $H$ for the corresponding vector
bundles: $\calg= \calo (G^\vee)$, $\calh=\calo(H^\vee)$. Together with
(\ref{obstr_theory_I}) we arrive at a morphism in the derived category
\[
  \ph^\bullet:\ [\calh\rightarrow \calg] \stackrel{\simeq}{\longrightarrow}
  (R\pi_*\Phi^*\calt_M)^\vee \longrightarrow \call_{Z/S}^\bullet\,.
\]
This is the local version of the free global normal space $\calf^\bullet
\rightarrow \call^\bullet_{\calc(M)/{\frak M}}$.  In fact, by
\cite[Prop.6.2]{behrendfantechi}
\begin{prop}
  $\ph^\bullet$ induces an isomorphism in $H^0$ and an epimorphism in
  $H^{-1}$.
  \qed
\end{prop}


\subsection{Analytic global normal cone}
We now review how to produce a cone (over $Z$) from a global normal space
$\phi^\bullet: [\calh\rightarrow\calg] \rightarrow \call^\bullet_{Z/S}$. It is
convenient to go over to linear fiber spaces over $Z$. In a hopefully
self-explanatory notation $\phi^\bullet$ thus becomes
$\Phi^\bullet: (L_{Z/S})_\bullet \rightarrow [G\rightarrow H]$.
Let $Z \hookrightarrow W$ be an embedding of $Z$ into a complex space $W$ that
is smooth over $S$. This is possible at least locally. Let $C_{Z|W}$ be the normal
cone of $Z$ in $W$, which is a closed subspace of the normal space $N_{Z|W}$,
the linear space associated to the conormal sheaf $\cali/\cali^2$, $\cali\subset
\calo_W$ the ideal defining $Z$. There is a natural morphism (in the derived sense)
of linear fiber spaces over $Z$ \cite[Cor.3.1.3]{illusie}
\[
  [T_{W|S}|_Z\rightarrow N_{Z|W}] \longrightarrow (L_{Z|S})_\bullet
\]
inducing isomorphism in $H^0$ and $H^1$. We are now in position to produce a
cone $C(\Phi^\bullet) \subset H$ by pushing forward the $T_{W|S}|_Z$-cone
$C_{Z|W} \subset N_{Z|W}$ under the composition of this morphism with
$\Phi^\bullet$ in the way defined in \cite[Ch.2]{si2}. $C(\Phi^\bullet)$ depends
only on the map in cohomology induced by $\Phi^\bullet$ \cite[Thm.3.3]{si2}.
\vspace{1ex}

From the free global normal space for unrigidified local models of $C(M)$
we thus get a well-defined closed analytic cone $C^H\subset H$.
To globalize first note that $G$ and $H$ have straightforward analogs
on $\calc(M)$. In fact, if $(\pi: \Gamma \rightarrow \calc(M),\tilde{\bf x})$ is the
universal (marked) curve with universal morphism
$\Phi: \Gamma \rightarrow M$ we put
\begin{eqnarray*}&&
   \call\ :=\ \omega_{\pi}(\underline{ x_1}+\ldots
  + \underline{ x_k}) \otimes {\Phi}^*
  \calh_M^{\otimes 3},\quad
  \caln\ :=\ {\pi}^*{\pi}_*({\Phi}^*\calt_M \otimes {
  \call}^{\otimes 3}) \otimes {\call}^{\otimes -3}\,,\\&&
  \calk\ :=\ \ker(\caln \rightarrow {\Phi}^*\calt_M)\,.
\end{eqnarray*}
Then $ G$, $ H$ are the orbibundles belonging to the orbisheaves
$R^1{ \pi}_* \calk$, $R^1{\pi}_* \caln$. 
To keep the notation within limits we stay with our previous symbols $G$, $H$,
$\Gamma$ etc. If we want to explicitely refer to uniformizing objects over
a local uniformizer $\hat U$ say, we will use the notation $H|_{\hat U}$ etc.

Now let $(C,{\bf x},\ph)\in \calc(M)$. A chart $\hat U$
for the analytic orbispace $\calc(M)$ was given by imposing incidence
conditions on an open $Z\subset \Hom_S(\calc,M)$, cf.\ the proof of
Proposition~\ref{CM_analyt_space}:
\[
  \hat U\ =\ (\ev_{k+1},\ldots, \ev_{k+l})^{-1} (H_1\times\ldots H_l)\,.
\]
Letting $\Phi^\bullet$ be associated to the universal objects on the unrigidified
chart $Z$ we define
\[
  {\hat C}^H\ :=\ C(\Phi^\bullet)|_{\hat U}\,,
\]
that is, the intersection of $C(\Phi^\bullet)$ with $H|_{\hat U}$ as closed
analytic subspaces of $H$ (meaning the bundle over $Z$). Since the
automorphism group of $(C, {\bf x},\ph)$ acts naturally on the whole
construction, ${\hat C}^H$ is $\Aut(C, {\bf x},\ph)$-invariant. With the
identification of $H|_{\hat U}$ with a local uniformizer of $H$ we may now
view ${\hat C}^H/ \Aut(C,{\bf x},\ph)$ as locally closed complex
subspace of $H$.
\begin{prop}
  The germ of ${\hat C}^H/ \Aut(C,{\bf x},\ph) \subset  H$ at any point of
  its support is independent of choices.
\end{prop}
\pf
To compare two local uniformizing systems $\hat U'\subset Z'\subset
\Hom_{S'}(\calc',M)$, $\hat U''\subset Z''\subset \Hom_{S''}(\calc'',M)$,
constructed from the insertion of $l'$ respectively $l''$ points
$y'_i$, $y''_i$ into $C$ (cf.\ the proof of Proposition~\ref {CM_analyt_space})
we may consider (appropriate shrinkings of) $Z'$ and $Z''$ as subspaces
of a common bigger space $Z''\subset \Hom_S (\calc,M)$. The latter
space is simply constructed by inserting the union of $\{y'_i\}$ and
$\{y''_i\}$, which is a tuple of $l\le l'+l''$ points. The point is that
for any increase of numbers of inserted points there is a forgetful map
(forgetting the additional points), and this induces morphisms
$Z\rightarrow Z'$, $Z\rightarrow Z''$. In fact, $Z$ is just a product of $Z'$
with an open set $D$ in $\cz^{l-l'}$ that rules the deformation of the
additionally inserted points. A similar statement is true for the universal
curve over $Z$, and this is obviously compatible with the evaluation
map. So the cone constructed from $Z$ is just a product of the cone
constructed from $Z'$ with $D$. Now the composition
\[
  Z'\ \longhookrightarrow Z\ \longrightarrow Z''
\]
induces the unique isomorphisms of the corresponding universal
deformations $U'$, $U$, $U''$ of $(C,{\bf x},\ph)$ and this shows the claim.
\qed

The locally closed subcones of (the underlying complex space of) $ H$
thus glue to a well-defined closed subcone $ C^H \subset H$ that is
locally pure dimensional of dimension $\rk G$ plus expected dimension.
\vspace{1ex}

The comparison with the symplectic treatment will happen on the
level of underlying spaces to which all relevant objects descend.
In the symplectic case this is the set of isomorphism classes of
stable holomorphic curves, whereas from the complex-analytic treatment
it also inherits the structure of a locally ringed space. In case $M$
is projective it follows either by construction or by using the universal
property that the latter space is exactly the analytic space corresponding to the
coarse moduli space underlying the DM-stack of stable curves in $M$
(which is in fact projective algebraic as shown in \cite{fultonpand}).

Now the (associated analytic) map from local \'etale covers of the
stack $C(M)_{\rm alg}$ to the coarse moduli space (\cite{behrendmanin},
after Prop.4.7) factorizes (locally analytically) through a smooth map to
our (even unrigidified) charts. Moreover, the GAGA-theorems from Chapter~5
show compatibility of the stack-theoretic global normal space and
our analytic global normal space $\ph^\bullet$. This proves
\begin{prop}
  The analytic cone $C\subset  H$ over the coarse moduli space
  associated to the stack-theoretic cone from \cite{behrend} coincides
  with the complex space underlying the complex orbispace $C^H$.
\qed
\end{prop}
To obtain the analytic analog of Behrend's virtual fundamental class $J(M)$
(an element in the Chow group $A_*(\calc(M))$\,) we finally need to intersect the
fundamental class $[ C^H]\in A_*( C^H)$ with the zero section
of $ H$. Algebraically this is done by applying a bivariant
class $\sigma^! \in A^*(\calc(M) \hookrightarrow  H)
\otimes_\gz\qz$. The existence of the latter in the category of
DM-stacks follows from the work of Vistoli. It results from the compatibility
of algebraic intersection theory with homology theory \cite[Ch.19]{fulton} that
the image of $J(M)$ in $H_*(\calc(M))$ (that we will also denote $J(M)$) is
nothing but
\begin{eqnarray}\label{alg_vfc}
  J(M)\ =\ [ C^H]\cap \Theta_{ H}\ \in\ H_*(\calc(M))\,,
\end{eqnarray}
where $[ C^H]$ is the fundamental class (of the underlying
coarse modulis spaces), and $\Theta_{ H} \in H_{\calc(M)}^{\rk H}( H, \qz)$
is the Thom class of the orbibundle $H$ (taking into account
multiplicities coming from the stabilizers of the local groups,
cf.\ \cite[\S1.2]{si1}).\footnote{
We could also define $J(M)$ as analytic Chow class, but since we want to
compare $J(M)$ with the homologically defined $\GW^M\in H_*(\calc(M))$ we
prefer to work in homology already at this point.}


\section{Analytic Kuranishi models}
A Kuranishi model for $\calc(M)$ is a locally closed embedding of a local
uniformizer of $\calc(M)$ into some $\cz^N$ won from a holomorphic Fredholm
map between complex Banach spaces having $\calc(M)$ as one of its fibers.
Finding Kuranishi models in our situation of an integrable complex structure on
$M$ is easier than in the general symplectic case, because we may restrict to
holomorphic maps on a large part of the curve, notably near the
singularities. This is due to the Stein property of open Riemann surfaces.
We thus begin with a study of moduli spaces of
holomorphic maps from an open Riemann surface.


\subsection{Spaces of holomorphic maps from open Riemann surfaces}
Throughout this section we fix an {\em open} Riemann surface $\Sigma$, whose
ideal boundary (denoted $\partial\Sigma$ by abuse of notation)
we assume to consist of circles only (no punctures). Denote  by
$\bar\Sigma= \Sigma\cup \partial\Sigma$ the corresponding Riemann
surface with boundary.
\vspace{1ex}

As a preparation let us generalize the decomposition
\[
  L_1^p(\Sigma;\cz)\ =\ L^p(\Sigma;\cz)\times \calo^{(c)}(\Sigma),
\]
that we used in \cite[\S4.3]{si1} for plane circular domains, to arbitrary
open Riemann surfaces. To that end let $K\in\Omega^1 (\Sigma \times\Sigma
\setminus\Delta)$ be a meromorphic 1-form with a simple pole along the diagonal
$\Delta$ extending to a neighbourhood of the boundary and without zeros.
$K$ exists by the Stein property of $\Sigma$.
Then the singular integral operator
\[
  T\gamma(z)\ :=\ \frac{1}{2\pi i} \lim_{\eps\rightarrow 0}
  \int_{\Sigma\setminus B_\eps(z)} K(z,w)\wedge \gamma(w)
\]
is a right-inverse for $\dbar: L_1^p(\Sigma;\cz) \rightarrow L^p(\Sigma;
\bar\Omega_\Sigma)$. Moreover, for any $f\in L_1^p(\Sigma;\cz)$ it holds
\[
  f\ =\ T\dbar f + Hf,\quad
  Hf(z)\ =\ \frac{1}{2\pi i}\int_{\partial\Sigma} K(z,w)\cdot f(w)\,.
\]
The proofs of these statements run as in the case $\Sigma=\Delta$. This
establishes the claimed decomposition of $L_1^p(\Sigma;\cz)$.
\vspace{1ex}

The rest of this section is devoted to the following result.
\begin{prop}\label{Stein_Hom}
  The space $\Hom^{(c)}(\Sigma,M)$ of holomorphic maps $\ph:\Sigma
  \rightarrow M$ that extend continuously to $\bar\Sigma$ has naturally a
  structure of complex Banach manifold. The evaluation map
  \[
    {\rm ev}:  \Hom^{(c)}(\Sigma,M)\times\Sigma \longrightarrow M,\quad
    (\ph,z)\longmapsto \ph(z)
  \]
  is holomorphic with respect to this complex structure. The tangent
  space at some $\ph\in\Hom^{(c)}(\Sigma,M)$ can naturally be identified
  with $\calo^{(c)}(\Sigma,\ph^*T_M)$, the space of holomorphic vector fields
  along $\Sigma$ extending continuously to $\bar\Sigma$.
  
  Moreover, with this complex structure, $\Hom^{(c)}(\Sigma,M)$ has
  the following ``universal property'': Let $T$ be a complex space and $\Phi:T\times
  \Sigma\rightarrow M$ be a holomorphic map extending continuously
  to $T\times \bar\Sigma$. Then there exists a unique holomorphic map
  $\rho: T\rightarrow \Hom^{(c)} (\Sigma,M)$ with
  \[
    \Phi\ =\ {\rm ev}\circ(\rho\times\id_\Sigma)\,.
  \]
  \vspace{-4ex}
  
\end{prop}
We put ``universal property'' in quotation marks because $\Hom^{(c)}(\Sigma,M)$
as an infinite dimensional space does not belong to the category of complex
spaces for which the property is tested (we would have to use Douady's
more involved notion of Banach analytic spaces to remedy this).

As a first step we observe that we can build up $\Sigma$
inductively starting from a disk by adding arbitrarily ``thin'' annuli or pairs
of pants. We are thus basically reduced to the case of arbitrarily thin 1-
or 2-connected plane domains by means of
\begin{lemma}
  Let $\Sigma= \Sigma_1\cup \Sigma_2$ with $\partial\Sigma_1\cap \partial \Sigma_2
  =\emptyset$ and $\ph\in \Hom^{(c)}(\Sigma,M)$. Put $\Sigma_{12}=
  \Sigma_1\cap\Sigma_2$ and $\ph_i=\ph|_{\Sigma_i}$, $\ph_{12}=\ph|_{\Sigma_{12}}$.
  Assume the proposition is true for $\Sigma_1$, $\Sigma_2$ and $\Sigma_{12}$
  locally around $\ph_1$, $\ph_2$ and $\ph_{12}$ respectively. Then the proposition
  is true for $\Sigma$ locally around $\ph$.
\end{lemma}
\pf
Let $\theta$ be an isomorphism of a neighbourhood of $\ph_{12}$ in $\Hom^{(c)}
(\Sigma_{12},M)$ with a neighbourhood of the origin in $\calo^{(c)}
(\Sigma_{12}; \ph_{12}^*T_M)$ whose differential at $\ph_{12}$ we assume to
be the identity. We consider the  map
\[
  \Xi: (\psi_1,\psi_2)\ \longmapsto\
  \theta(\psi_1|_{\Sigma_{12}})-\theta(\psi_2|_{\Sigma_{12}})
\]
from a neighbourhood of $(\ph_1,\ph_2)$ in $\Hom^{(c)}(\Sigma_1,M)
\times\Hom^{(c)}(\Sigma_2,M)$ to $\calo^{(c)} (\Sigma_{12}; \ph_{12}^*T_M)$.
The derivative of this map is nothing but the \v Cech-differential on sections
of $\ph^*T_M$ associated to the covering $\Sigma= \Sigma_1\cup \Sigma_2$.
We claim that the latter is a split submersion. To this end let $T^1$ and $T^2$
be right-inverses to $\dbar$ on $\ph_1^*T_M$ and $\ph_2^*T_M$ as above.
For simplicity we assume the meromorphic 1-forms (on $\Sigma_i\times\Sigma_i$)
involved in the definition to be restrictions of a meromorphic 1-form $K$
on $\Sigma\times\Sigma$. Note that the pull-backs of $T_M$
are trivial as holomorphic vector bundles, so the definition
of a right-inverse to $\dbar$ extends to this case by use of a
vector valued analog with respect to a holomorphic trivialization of $\ph^*T_M$
(which we may take to extend continuously to $\bar\Sigma$). Let $1\equiv\rho_1
+\rho_2$ be a partition of unity subordinate to $\Sigma_i$. A splitting
of the \v Cech-complex is given by
\begin{eqnarray*}
  \calo^{(c)}(\Sigma_{12};\ph_{12}^*T_M)& \longrightarrow&
  \calo^{(c)}(\Sigma_1;\ph_1^*T_M)\times \calo^{(c)}(\Sigma_1;\ph_1^*T_M)\\
  v&\longmapsto& (\rho_2\cdot v-T^1(v\dbar\rho_2),
  -\rho_1\cdot v+T^2(v\dbar\rho_1))\,.
\end{eqnarray*}
In fact, letting $T$ be the singular integral operator on $\Sigma$ belonging
to $K$, we obtain
\[
  \rho_2\cdot v-T^1(v\dbar\rho_2)+\rho_1\cdot v-T^2(v\dbar\rho_1))\ =\ 
  v-T(v\dbar (\rho_1+\rho_2))\ =\ v\,.
\]
In view of the holomorphic implicit function theorem
this shows the statement about the structure of complex Banach manifold
and the description of its tangent spaces. The ``universal property'' is obvious
because for any complex space $T$ a holomorphic map
$T\times\Sigma\rightarrow M$ is equivalent to the giving of a pair of
holomorphic maps $T\times\Sigma_i\rightarrow M$ coinciding on
$T\times \Sigma_{12}$.
\qed
\vspace{1ex}

Next we treat the case that the map $\ph$ is ``small'' compared to the modulus
of $\Sigma$, as made precise in the following lemma.
\begin{lemma}\label{thindomains}
  Let $\ph\in \Hom^{(c)}(\Sigma,M)$ and assume there is a covering $U_i$
  of $\Sigma$ by finitely many open sets having piecewise smooth boundary
  in $\bar\Sigma$ and such that (1)~there exists a smooth partition of unity
  $\{\rho_i\}$ subordinate to $U_i$
  (2)~$U_i\cap U_j\cap U_k=\emptyset$ for any
  three pairwise different indices $i,j,k$ (3)~for any $i$ the closure of $\ph(U_i)$
  is contained in a holomorphic coordinate chart of $M$. Then the assertions of the
  proposition hold locally around $\ph$.
\end{lemma}
Assumption (1) means that a shrinking of $\{U_i\}$ is still a covering of $\Sigma$.
This requires the closures of $U_i$ in $\bar\Sigma$ to meet in one-dimensional
subsets on $\partial \Sigma$.\\[1ex]
\pf
Let $\gamma_i:M \supset
W_i\rightarrow\cz^n$ denote holomorphic coordinates on $M$
with $\cl\ph(U_i)\subset W_i$. We consider the gluing map
\[
  \Xi: \prod_i \Hom^{(c)}(U_i,W_i)\longrightarrow
  \prod_{i<j} \Hom^{(c)}(U_i\cap U_j,\cz^n)\,,\quad
  (\psi_i)_i\longmapsto \Big(\gamma_i\circ\psi_i|_{U_j}
  -\gamma_i\circ\psi_j|_{U_i}\Big)_{ij}\,.
\]
Clearly, the spaces $\Hom^{(c)}(U_i,W_i)$ are open sets in complex Banach spaces
via $\gamma_i$. As in the previous lemma we want to show that the differential
of $\Xi$ at $(\ph|_{U_i})_i$ is a split submersion. This differential is the map
\[
  \prod_i \calo^{(c)}(U_i;\ph^*T_M)\longrightarrow
  \prod_{i<j} \calo^{(c)}(U_i\cap U_j,;\cz^n)\,,\quad
  (v_i)_i \longmapsto (D\gamma_i)(v_i-v_j)\,.
\]
Note that the differential of $\gamma_i$ is an isomorphism for any $i$.
The requested splitting can thus be defined by
\[
  (D\gamma_i (v_{ij}))_{ij}\longmapsto
  \Big(\sum_j \rho_j v_{ij} - T^i(v_{ij}\dbar \rho_j)\Big)_i\,.
\]
As in the proof of the previous lemma we use a holomorphic trivialization
of $\ph^*T_M$ to construct a right-inverse $T$ to $\dbar$, the restriction of
which to sections with support in $U_i$ is $T^i$. The verification that this is
indeed a splitting of the differential and the ``universal property'' run also
as above. We conclude that $\Xi^{-1}(0)$ is a chart for $\Hom^{(c)}(\Sigma,M)$
at $\ph$ with tangent space $\calo^{(c)}(\Sigma, \psi^*T_M)$ at $\psi$.
The universal property is again obvious.
\qed
\vspace{1ex}

\noindent
{\em Proof of Proposition~\ref{Stein_Hom}.}\ 
Locally around some $\ph\in \Hom^{(c)} (\Sigma,M)$ we construct the
complex Banach manifold structure in the way already indicated: Write
$\Sigma$ as union of open Riemann surfaces $\Sigma_i$ with the following
properties: (1) Each $\Sigma_i$ is either a disk, an annulus, or a pair of
pants (a two-connected plane domain)
(2) For any $i\neq j$ the intersection $\Sigma_i\cap \Sigma_j$ is
a union of (one or two) annuli (3) Each $\Sigma_i$ has a covering by open
sets as required by Lemma~\ref{thindomains} (the same holds then also
true for the pairwise intersections). Applying the two lemmata we obtain
the statement of the proposition locally around $\ph$.

It remains to show
independence of choices (in particular of the decomposition $\Sigma=\bigcup
\Sigma_i$) and biholomorphicity of changes of coordinates.
So let $V=\Xi^{-1}(0)$, $V'=(\Xi')^{-1}(0)$ be two such charts with $V\cap V'
\neq\emptyset$ as subsets of $\Hom^{(c)}(\Sigma,M)$. From the holomorphic
evaluation map $V'\times \Sigma\rightarrow M$, upon restriction to $\Sigma_i$,
we obtain maps
\[
  V'\ \longrightarrow\ \Hom^{(c)}(\Sigma_i,M)\,.
\]
These are holomorphic maps of Banach manifolds and compatible
on intersections $\Sigma_i\cap\Sigma_j$. We thus obtain a holomorphic
map of Banach manifolds $V'\rightarrow V$. The same reasoning with
$V,V'$ exchanged shows that this map is invertible.
\qed


\subsection{Localization of variations of complex structure}
As a further preliminary to the construction of analytic Kuranishi
models we want to show that any deformation $q:X\rightarrow S$ of a closed
Riemann surface $\Sigma=X_0$ is obtained by changing the complex structure
on {\em arbitrarily small} open sets. To see this take a meromorphic function $f$ on
(a neighbourhood of the central fiber of) $X$ exhibiting $\Sigma$ as a
branched covering of $\pr^1$. Then nearby fibers are also branched
coverings of $\pr^1$, and the set of branch points varies
holomorphically with $s$. Let $V_\eps\subset \pr^1$ be an
$\eps$-neighbourhood (in any metric on $\pr^1$) of the
branch locus of $\Sigma\rightarrow \pr^1$. For sufficiently small $s$ the branch
locus of $X_s\rightarrow \pr^1$ is still contained in $V_\eps$. For any such $s$
the map $f$ induces an isomorphism $X_s\setminus\ph^{-1} (V_\eps) \simeq
\Sigma\setminus\ph^{-1} (V_\eps)$. Note that if we are only interested in the
germ of the deformation we may take $\eps$ arbitrarily small. For any branch point
$P\in\Sigma$ let $z$ be a local holomorphic function defined in a neighbourhood
of $P$ in $X$ that restrict to a local holomorphic coordinate on $\Sigma$.
Provided $\eps$ is sufficiently small $z$ will be defined on the closure of the
connected component of $f^{-1}(V_\eps)\cap \Sigma$ containing $P$. So $z$
will map this connected component to a simply connected plane domain, and
the same holds for $f^{-1}(V_\eps)\cap X_s$, $s$ sufficiently small. We thus
obtain a description of $X\rightarrow S$ as gluing of an open Riemann surface
$\Sigma\setminus f^{-1}(\cl V_{\eps/2})$ with a family of plane domains $W_j(s)$
of the form $z(f^{-1}(V_\eps))$. The gluing is along the union of constant
annuli $\Sigma\cap f^{-1}(V_\eps\setminus \cl V_{\eps/2})$.

The point is that
the plane domains vary holomorphically with $s$ in the following sense: Assume
$P$ is a branch point of order $b$. Let $g: \Delta\rightarrow \pr^1$ be the
composition of $\Delta\rightarrow\Delta, t\rightarrow t^b$ with an isomorphism
$\Delta\rightarrow V_\eps$ (extending holomorphically to the closure) and
mapping $0$ to $f(P)$. Then $g$ factorizes over $f$ via a holomorphic map
$h:\Delta\rightarrow\Sigma$ with $h(0)=P$. Note that $h$ is a local uniformizer for
$\Sigma$ at $P$ extending holomorphically to $\cl\Delta$ and inducing a
diffeomorphism of $S^1$ with the boundary of the connected component of
$f^{-1}(V_\eps)$ containing $P$. So $z\circ h$ provides an isomorphism of $\Delta$
with $W_j(0)$. Going over to nonzero $s$ the same reasoning applies but we
will have to restrict $g$ (thus $h$) to an annulus $A$ with the property that
$g(A)$ does not meet the branch locus of $f:X_s\rightarrow \pr^1$. We thus see
that the boundary of a family $W_j(s)$ is the image of a map $h:S^1\times S
\rightarrow \cz$ with $h_\ph=h|_ {\{\ph\}\times S}$ {\em holomorphic}. We will
refer to this property by saying that $W_j(s)$ is a {\em holomorphic family of
simply connected plane domains}.

Taking into account also deformations of nodes ($Z_t=\{ (z,w)\in\Delta \times
\Delta \mid zw=t\}$) we obtain:
\begin{lemma}\label{decomposition}
  Let $p: X\rightarrow S$ be a deformation of a prestable curve $C$.
  Then, possibly after shrinking $S$ there exists a decomposition
  $X=U_0\cup \bigcup_{i=1}^l U_i \cup \bigcup_{j=1}^m W_j$ with
  \begin{itemize}
  \item
    $U_0=S\times\Sigma$, where $\Sigma$ is a non-compact Riemann surface
  \item
    $U_i(s)=Z_{t_i(s)}$ for some $t_i\in\calo(S)$ with $t_i(0)=0$
  \item
    $W_j\rightarrow S$ is a holomorphic family of simply connected plane domains
  \item
    $U_i\cap U_{i'}=U_i\cap W_j= W_j\cap W_{j'}=\emptyset$ for any $i,i',j,j' \ge 1$,
    $i\neq i'$, $j\neq j'$; $U_0\cap U_i\simeq S\times (A_i\cup A'_i)$;
    $U_0\cap W_j\simeq S\times A''_j$ for some annuli $A_i, A'_i,A''_j$.
  \end{itemize}
  Moreover, the $U_i$ ($i>0$) and $W_j$ can be chosen arbitrarily small.
\qed
\end{lemma}
\vspace{1ex}

The fact that the boundary of $W_i(s)$ varies holomorphically can be used
to give a holomorphic trivialization of $\coprod_{s\in S} \calo^{(c)}(W_i(s))$
via Cauchy integrals. Recall that for any simply connected domain $G\subset \cz$
with piecewise smooth boundary the Cauchy integral $F(z)=\frac{1}{2\pi i}
\int_{\partial G} \frac{f(w)}{w-z}dw$ induces a decomposition
of $C(\partial G)$ into the direct sum of $\calo^{(c)}(G)$
and the complex conjugate (modulo constants) $\overline{\calo^{(c)}(G)}/\cz$.
Choosing a biholomorphism $G\simeq\Delta$ (this extends to a
homeomorphism of the boundaries) the decomposition is exhibited
by Fourier series with only nonnegative and such with only negative
Fourier components: $C(S^1)=C_{\ge 0}(S^1)\oplus C_{<0} (S^1)$.
Let $h:S^1\times S\rightarrow \cz$ describe the holomorphic variation of
$W(s)=W_i(s)$, $h_s=h|_{S^1\times\{s\}}$. Note that $h|_{S^1\times\{0\}}$
extends to a biholomorphism $\Delta\simeq W(0)$.
\begin{prop}\label{bundlestructure}
  Possibly after shrinking $S$ the map
  \begin{eqnarray*}
    &S\times\calo^{(c)}(W(0))\simeq S\times C_{\ge 0}(S^1)\ \longrightarrow\ 
    \coprod_{s\in S} \calo^{(c)}(W_j(s))\,,\quad
    (s,f)\ \longmapsto\ (s,\Phi_{f,s})\\[2pt]
    &\displaystyle\Phi_{f,s}(z)=\frac{1}{2\pi i}\int_{\partial W(s)}
    \frac{f(h_s^{-1}(w))}{w-z} dw\,,\quad z\in W(s)
  \end{eqnarray*}
  is bijective. Moreover, as a function of $(s,z)$, $\Phi_{f,s}(z)$ is holomorphic.
\end{prop}
\pf
For each $s$ write $\Phi_s(f)=\Phi_{f,s}$. Then $\Pi_s:= h_s^ *\circ \Phi_s$ is
a continuous family of projections on $C(S^1)$. The image of $\Pi_s$ corresponds
to the space of boundary values of holomorphic functions on $W(s)$, via
pull-back by $h_s$. But at $s=0$ this is $C_{\ge 0}(S^1)$ and so $\Pi_s$
induces an isomorphism of $C_{\ge 0}(S^1)$ with $\img\Pi_s$ for sufficiently
small $s$. This proves the claimed bijectivity.

Given any $f\in C(S^1)$ the holomorphicity of $\Phi_{f,s}(z)$ in
$(s,z)$ is obvious by writing
\[
  \Phi_{f,s}(z)=\frac{1}{2\pi}\int_0^{2\pi}\frac{f(e^{i\ph})}
  {h_s(e^{i\ph})-z}\, e^{i\ph}\,d\ph\,.
\]
\vspace{-4ex} 

\qed


\subsection{Construction of analytic Kuranishi models}
We now want to make the construction of (local) Kuranishi models in \cite{si1}
holomorphic. At $(C, {\bf x}, \ph)\in \calc(M)$ the construction worked as
follows: Let
\[
  (q:\calc\rightarrow S, \underline{\bf x}:S\rightarrow
  \calc\times_S\ldots \times_S \calc)
\]
be an analytically semiuniversal deformation of $(C,{\bf x})$, $C=q^{-1}(0)$. We
showed in op.\ cit.\ that the space of Sobolev maps (with weights at the
nodes) from the $C_s= q^{-1}(s)$ to $M$
\[
  {\check L}_1^p (\calc/S;M)\ =\ \coprod_{s\in S} {\check L}_1^p (C_s;M)
\]
is a Banach manifold, locally isomorphic to $S\times V$, $V\subset {\check
L}_1^p (C; \ph^*T_M)$ open. The map $\ph\mapsto \dbar \ph$ is locally a family
of differentiable Fredholm maps
\[
  H_s=H|_{\{s\}\times V},\quad
  H: S\times V\longrightarrow E:= {\check L}^p(C;\ph^* T_M
  \otimes \bar\Omega)\,.
\]
Choose a linear projection $\prj_Q: E\rightarrow Q$ with finite-dimensional
kernel and with $\prj_Q \circ D H_0$ surjective. The implicit function theorem
with parameter $s$ then shows that $\tilde Z=(\prj_Q \circ H)^{-1} (0)$ is a
finite-dimensional manifold near $(0,0)\in S\times V$. Provided $(C,{\bf x})$
is stable, a local uniformizer of $\calc(M)$ is locally given as zero locus of
$\kappa|_{\tilde Z}$, $\kappa=H-\prj_Q \circ H =\prj_F
\circ H$, $\prj_F:E \rightarrow F$ the projection with kernel $Q$ onto a finite
dimensional subspace $F\subset E$ with $E= F+ \img D H_0$. So $F$ spans the
cokernel of $DH_0$. If $(C,{\bf x})$ is not stable one finally has to choose a slice
to the action of $\Aut^0 (C,{\bf x})$ on $S\times V$.
\vspace{1ex}

Instead of making this construction holomorphic in each step we restrict to a
set of maps in ${\check L}_1^p (\calc/S; M)$ that are already largely holomorphic.
We take a decomposition of $C$ as provided by Lemma~\ref{decomposition} applied
to $q:\calc\rightarrow S$. Without loss of generality we may assume that
each irreducible component of $C$ contains at least one $W_j$ that
remains constant under the deformation, say $j\in\{ 1,\ldots, m'\}$.
We assume that the closure of the image of $U_i$ ($i>0$)
and of $W_j$ under $\ph$ map into holomorphic coordinate charts
\[
  \gamma_i: M\supset M_i\rightarrow\cz^n\,,\quad i=1,\ldots,l+m\,.
\]
So $M_i$ are open sets in $M$ containing $\cl\ph(U_i(0))$ respectively $\cl
\ph(W_{i-l}(0))$ (for $i>l$). Our ambient Banach manifold will be the subset
of ${\check L}_1^p (\calc/S;M)$ of maps that are holomorphic on
$\bigcup_{i\ge0} U_i\cup \bigcup_{j=m'+1}^m W_j$. So the flexibility
provided by $L_1^p$-maps survives only on the part of $W_1,\ldots,
W_{m'}$ not meeting $\Sigma=U_0$. The space of such maps
can be described as fiber over $0$ of a holomorphic map $\Xi:B\rightarrow B'$
of complex Banach manifolds. The domain $B$ is the fibered product over $S$
of the spaces
\begin{eqnarray*}
  &S\times\Hom^{(c)}(\Sigma,M)\,;\quad \Hom_S^{(c)}(U_i, M_i)\,,i=1,\ldots,l\,;\\
  &S\times L_1^p (W_j(0); M_{l+j})\,,j=1,\ldots,m'\,;\quad
  \Hom_S^{(c)}(W_j, M_{l+j})\,,j=m'+1,\ldots,m\,,
\end{eqnarray*}
(we use Proposition~\ref{Stein_Hom}) and $B'$ will be
\[
  \prod_{i=1}^l \Hom^{(c)}(A_i\cup A'_i,\cz^n)
  \times\prod_{j=1}^{m'} L_1^p (A''_j; \cz^n)
  \times\prod_{j=m'+1}^m \Hom^{(c)}(A''_j,\cz^n)\,.
\]
Here we wrote $\Hom_S^{(c)}(U_i, M_i)$ for $\coprod_{s\in S}
\Hom^{(c)}(U_i(s), M_i)$ that in view of the identification of holomorphic function
spaces on deformation spaces of nodes given in \cite[\S4.2]{si1}
is an open set in the product of $S$ with a complex Banach space
via $\gamma_i$. By Proposition~\ref{bundlestructure} a similar remark applies to
$\Hom_S^{(c)}(W_j, M_{l+j}) =\coprod_{s\in S} \Hom^{(c)} (W_j(s), M_{l+j})$. 
The map $\Xi$ is the obvious nonlinear version of the first \v Cech differential.
It sends $(\psi_0,\psi_i, \psi'_j)$ to
\[
  \Big(\gamma_i\circ\psi_i-\gamma_i\circ\psi_0,
  \gamma_j\circ\psi'_j-\gamma_j\circ\psi_0\Big)\,,
\]
and this is clearly a holomorphic map.
\begin{prop}
  In a neighbourhood of $(C,{\bf x},\ph)$,
  $\calb:= \Xi^{-1}(0)$ is a complex Banach manifold lying smoothly over $S$.
\end{prop}
\pf
In view of the implicit function theorem we have to show that the differential
$D \Xi_0$ is a split submersion, where $\Xi_0$ is the restriction of
$\Xi$ to $s=0$. With the natural identification of the tangent
space relative $S$ of $B$ at $\ph$ with
\[
  \calo^{(c)}(\Sigma;\ph^*T_M)\times \prod_{i=1}^l \calo^{(c)}(U_i, \ph^*T_M)
  \times\prod_{j=1}^{m'} L_1^p (W_j(0); \ph^*T_M)
  \times\prod_{j=m'+1}^m \calo^{(c)}(W_j(0), \ph^*T_M)
\]
the differential is isomorphic to the linear \v Cech differential, mapping
$(v_0,v_i,v'_j)$ to $(v_i-v_0,v'_j-v_0)$, as element in
\[
  \prod_{i=1}^l \calo^{(c)}(A_i\cup A'_i;\ph^*T_M)
  \times\prod_{j=1}^{m'} L_1^p (A''_j; \ph^*T_M)
  \times\prod_{j=m'+1}^m \calo^{(c)} (A''_j, \ph^*T_M)\,.
\]
To define the splitting let $T$ be a right-inverse to the $\dbar$-operator
on $\ph^*T_M$ restricted to $C^0= C^0_0$, where
\[
  C^0_s:= \bigcup_{i=0}^l U_i(s)\cup \bigcup_{j=m'+1}^m W_j(s)\,,
\]
constructed as singular
integral operator via a holomorphic trivialization as in the last section.
Let $T^i$ be the restriction of $T$ to $U_i(0)$ (for $i=0,\ldots,l$)
and to $W_{m'+i-l}(0)$ respectively (for $i=l+1,\ldots, l+m-m'$). For brevity
we put $H^i=\id-T^i\circ \dbar$. For $j=1,\ldots,m'$ we also need a
(complex linear) extension map
\[
  \Omega^j: L_1^p (U_0(0)\cap W_j(0); \ph^*T_M)\ \longrightarrow\
  L_1^p (W_j(0); \ph^*T_M)\,,
\]
that is, a right-inverse to the corresponding
restriction map. The existence of $\Omega^j$ is a standard
fact of Sobolev theory. Finally let $\rho_i$, $\rho'_j$ be a partition of unity
on $C$ subordinate to our covering. A right inverse to $D\Xi_0$ is now
given by sending $(w_i,w'_j)$ to $(v_0,v_i,v'_j)$ with
\[\begin{array}{rcll}
  v_0&=& -H^0 (\sum_{i=1}^l \rho_i w_i+\sum_{j=m'+1}^m \rho_j w'_j) \\
  v_i&=& H^i(\rho_0 w_i)& i=1,\ldots,l\\
  v'_j&=& \Omega^j(w'_j+(v_0|_{A''_j}))& j=1,\ldots,m'\\
  v'_j&=& H^{l+j-m'}(\rho'_0 w'_j)& j=m'+1,\ldots,m\,.
\end{array}\]
The verification that this is indeed a splitting is straightforward (as in the
proofs of the lemmata in the previous section).
An application of the implicit function theorem completes the proof.
\qed

In the sequel we identify $(\psi_0, \psi_i,\psi'_j)\in \calb$ with the induced
$L_1^p$-map $\psi: C_s\rightarrow M$ that is holomorphic on $C^0_s$.
Choosing a biholomorphism $z_j:\Delta^j\simeq\Delta$ for
the complement $\Delta^j$ of
$\cl A''_j$ in $W_j(0)$ ($j=1,\ldots,m'$), the (non-linear) $\dbar$-equation
$\ph \mapsto \dbar\ph$ can now be viewed as holomorphic map
\[
  \Theta: \calb \longrightarrow \prod_{j=1}^{m'} L^p(\Delta^j; \cz^n), \quad
  (\psi_0,\psi_i,\psi'_j)\longmapsto
  \Big(\frac{\partial}{\partial \bar z_j}\,\gamma_j\circ\psi'_j \Big)_{j=1,\ldots,m'}\,.
\]
\begin{prop}\label{Tbar_quasiiso}
  There is a natural map from
  the differential at $(C,{\bf x},\ph)$ of $\Theta$ relative $S$ to
  \[
    \dbar:\ {\check L}_1^p(C;\ph^* T_M) \longrightarrow
    {\check L}^p(C;\ph^*T_M\otimes\bar\Omega)
  \]
  inducing isomorphisms of kernels and cokernels.
  In particular $\Theta$ is Fredholm at $(C,{\bf x},\ph)$.
\end{prop}
\pf
Let us write $\Theta_0$ for the restriction of $\Theta$ to the central fiber $\calb_0$
of $\calb$ over $0\in S$. From the proof of the last proposition the tangent space
of $\calb_0$ at $\psi: C_s\rightarrow M$ can be identified with
\[
  V:=\ \{v\in{\check L}_1^p (C;\ph^*T_M)\,|\, v|_{C^0}\in\calo(C^0;\ph^*T_M)\}\,.
\]
Thus
\[
  D\Theta_0:\ V\ \longrightarrow\ \prod_{j\le m'} L^p(\Delta^j;\cz^n)\,,\quad
  v\ \longmapsto\ \Big(D\gamma_j(\frac{\partial v}{\partial \bar z_j}) \Big)_j\,.
\]
Multiplying the components on the right-hand side by $d\bar{z_j}$ and
extending trivially by zero as section of $\ph^*T_M \otimes\bar\Omega$ we
obtain a commutative square
\begin{eqnarray}\label{comm1}
  &\begin{array}{ccc}
  V&\stackrel{D\Theta_0}{\llongrightarrow}&
  \prod_{j\le m'} L^p(\Delta^j;\cz^n)\\[2pt]
  \longdownarrow&&\longdownarrow\\[6pt]
  {\check L}_1^p (C;\ph^*T_M)&\stackrel{\dbar}{\llongrightarrow}&
  {\check L}^p(C;\ph^*T_M\otimes \bar\Omega)
  \end{array}&
\end{eqnarray}
that we claim to induce isomorphisms of kernels and cokernels of the
horizontal arrows. The diagram certainly gives rise to
isomorphisms on the kernels, both being equal to $H^0(C; \ph^*T_M)$.
To prove injectivity on cokernels let $\gamma\in \prod_{j\le m'}
L^p(\Delta^j,\cz^n)$ and assume that its trivial extension $\tilde\gamma$ can be
written as $\dbar v$ with $v\in \check L_1^p(C; \ph^*T_M)$. Then $\dbar v=0$ away
from $\bigcup \Delta^j$ and so $v\in V$.

As for surjectivity, the soft resolution of $\calo(\ph^* T_M|_{C^0})$ by sheaves
of Sobolev sections (restriction of \cite[Prop.4.5]{si1} to $C^0$)
\[
  0\longrightarrow \calo(\ph^*T_M|_{C^0}) \longrightarrow
  {\check \call}_1^p(\ph^*T_M|_{C^0})\stackrel{\dbar}{\longrightarrow}
  {\check\call}^p(\ph^*T_M|_{C^0})\longrightarrow0
\]
together with $H^1(C^0;\ph^*T_M)=0$ ($C^0$ is Stein!) shows the
surjectivity of
\[
  \dbar:\ {\check L}_1^p(C^0;\ph^*T_M)\longrightarrow
  {\check L}^p(C^0;\ph^*T_M\otimes\bar\Omega)\,.
\]
So for any $\gamma\in {\check L}^p(C^0;\ph^*T_M\otimes\bar\Omega)$
there is a solution to the equation
\[
  \dbar v_0\ =\ \gamma|_{C^0}\,.
\]
Let $v\in {\check L}_1^p (C;\ph^*T_M)$ be an extension of $v_0$ to all of
$C$. Then $\gamma-\dbar v$ has support in $\bigcup_{j\le m'} \Delta^j$, hence
is in the image of the right-hand vertical map of (\ref{comm1}). This shows
surjectivity of (\ref{comm1}) on cokernels.
\qed

By the proposition, $\Theta^{-1}(0)$ is thus given as the fiber of a holomorphic
Fredholm map between complex Banach manifolds, hence has naturally the
structure of a (finite dimensional) complex space.
\begin{theorem}\label{hol_str}
  The germs of $\Theta^{-1}(0)$ and of $\Hom_S (\calc;M)$ at $(C,{\bf x},\ph)$
  are canonically isomorphic.
\end{theorem}
\pf
We have to check the universal property of the hom-functor for $\Theta^{-1}(0)$.
So let $T\rightarrow S$ be a morphism of complex spaces mapping a distinguished
point $0\in T$ to $0\in S$ and write $\calc_T =T\times_S \calc$.
To any morphism $\Phi: \calc_T\rightarrow M$ inducing $\ph$ on the central fiber
we claim the existence of a unique (germ of) morphism $\Lambda:
T\rightarrow \Theta^{-1}(0)$ such that $\Phi$ factorizes over the evaluation
morphism $\Theta^{-1}(0)\times_S \calc \rightarrow M$ via
$\Lambda\times \id_\calc$.

First we observe that the morphism $\Phi$ from $\calc_T$ to $M$ is equivalent
to the giving of a tuple of $S$-morphisms from $T$ to the following spaces:
$S\times\Hom^{(c)}(\Sigma,M)$, $\Hom_S^{(c)}(U_i,M_i)$ ($i=1,\ldots,l$),
$\Hom_S^{(c)}(W_j,M_{l+j})$ ($j=1,\ldots,m$), such that the composition
with the analog $\Xi_h$ of $\Xi$ from the product $B_h$ (``$h$'' for
``holomorphic'') of these spaces to
\[
  B'_h:= \prod_{i=1}^l \Hom^{(c)}(A_i\cup A'_i,\cz^n)
  \times\prod_{j=1}^m \Hom^{(c)}(A''_j,\cz^n)
\]
is the zero morphism. Recall that the term ``$S$-morphism'' means compatibility
with the morphisms to $S$ that the relevant spaces do possess (including $T$).
Now consider the following diagram of complex Banach manifolds
\[\begin{array}{ccccl}
  &&B_h&\stackrel{\Xi_h}{\longrightarrow}& B'_h\\[2pt]
  &&\longdownarrow&&\ \ \longdownarrow\\[6pt]
  \calb&\longrightarrow&B&\stackrel{\Xi}{\longrightarrow}& B'\\[2pt]
  \diagl{\Theta}&&\diagl{\dbar}&&\ \ \diagr{\dbar}\\[6pt]
  \prod_{j=1}^{m'} L^p(\Delta^j;\cz^n)&\longrightarrow&
  \prod_{j=1}^{m'}  L^p(W_j(0);\cz^n)&\longrightarrow&
  \prod_{j=1}^{m'} L^p(A''_j;\cz^n)
\end{array}\]
The right-hand horizontal maps are all split submersions, the left-hand horizontal
and upper vertical maps are closed embeddings (the kernels of the following maps)
and the lower vertical maps are given by $\dbar$-operators. What we have
just said about families amounts to saying that the fiber over $0$ (as ringed space)
of the non-linear Fredholm map $B_h\rightarrow B'_h$
represents the Hom-functor under consideration. Now
a morphism to this fiber is equivalent to a morphism
to $B$ which composed with each of the maps to $B'$ and to
$ \prod  L^p(W_j(0);\cz^n)$ is the zero morphism. So this is nothing but a
morphism to $\calb$ whose composition with $\Theta$ to
$\prod L^p(\Delta^j;\cz^n)$ is zero, which in turn is the same as a morphism to
the fiber of $\Theta$ over $0$. All these morphisms are to be understood in the
category of ringed spaces. This shows the universal property of the hom-functor
for $\Theta^{-1}(0)$.
\qed

We do not bother about explicit local analytic Kuranishi models or analytic
rigidification here, because we will never need this. The former can be easily
established by the usual procedure of choosing a finite-dimensional
complex linear subspace in $\prod_{j\le m'} L^p(\Delta^j;\cz^n)$
spanning $\coker D\Theta_0$; a holomorphic rigidification on the other hand
can be provided by requiring incidence of additional marked points with local
analytic divisors in $M$ as in Proposition~\ref{CM_analyt_space}. Instead we will
show in Section~4.1 how to incorporate the local construction given here in the
global construction of \cite{si1}.


\section{Limit cones}
\subsection{Cone classes}
This section deals with a purely topological question in finite dimensions. Let
$N$ be a finite dimensional, oriented, topological orbifold, $q: F\rightarrow N$
an orbibundle and $s$ a continuous section of $F$.\footnote{
More generally, we may take $N$ to be an oriented, $\qz$-homology manifold
and $F\rightarrow N$ a not necessarily locally trivial cone over $S$, i.e.\ $F$
have a continuous, fiber-preserving action of the multiplicative semigroup
$\rz_{\ge0}$ that is proper away from the zero section $0\cdot F\simeq N$}
For any $l >0$ let $p_l: N\times \rz^l \rightarrow N$ be the projection and let
$s_l\in \Gamma(N\times(\rz^l \setminus\{0\}), p_l^*F)$ be defined by
\[
  s_l(x,v)\ =\ |v|^{-1}\cdot s(x)\,.
\]
For any $v\in \rz^l\setminus\{0\}$ the restriction of the graph $\Gamma_{s_l}$
of $s_l$ to $N\times \{v\}$ is $\Gamma_{|v|^{-1} \cdot s}$.
\begin{defi}\rm
Let $\overline{\Gamma_{s_l}}$ be the closure of $\Gamma_{s_l}$ in $p_l^* F$.
The {\em cone $C(s)\subset F$ associated to $s$} is defined by
\[
  C(s)\ :=\ \overline{\Gamma_{s_l}}\cap (F\times\{0\})\,.
\]
\vspace{-8ex}

\qed
\end{defi}
$C(s)$ is independent of $l$ and lies over the zero locus of $s$:
\begin{prop}\label{cones_prop}
1)\ \ $C(s)=\lim_{t\rightarrow\infty} \Gamma_{t\cdot s}$ in the sense of
Hausdorff convergence of closed sets (hence $C(s)$ is independent of
$l$)\\[1ex]
2)\ \ $q(C(s))=Z(s)$.
\end{prop}
\pf
1)\ Recall that convergence $\Gamma_{t\cdot s} \stackrel{t\rightarrow\infty}
{\longrightarrow} C(s)$ means two things:
\begin{itemize}
\item
  $C(s)=\bigcap_{t_0}\cl(\bigcup_{t\ge t_0}\Gamma_{t\cdot s})$
\item
  For any compactum $K\subset F$ and any neighbourhood $U$ of $C(s)\subset
  F$ there exists a $t_0$ with
  \[
    \Gamma_{t\cdot s}\cap K\subset U\quad\mbox{for any}\quad t>t_0.
  \]
\end{itemize}
Both properties follow easily from the corresponding facts for $\overline{
\Gamma_l}$: Since for any $t_0>0$ the restriction of the projection $P_l:
p_l^*F = F\times\rz^l \rightarrow F$ to $F\times \bar B_{t_0^{-1}}(0)$
is proper and $F$ is locally compact we have
\begin{eqnarray*}
  \cl_F\Big(\bigcup_{t\ge t_0}\Gamma_{t\cdot s}\Big)&=&
  P_l\Big[\cl_{F\times\rz^l}(\Gamma_{s_l}\cap
  (F\times\bar B_{t_0^{-1}}(0))\Big]\\
  &=& P_l\Big[\overline{\Gamma_{s_l}}\cap
  (F\times\bar B_{t_0^{-1}}(0))\Big]\,,
\end{eqnarray*}
where we indicated with subscripts in which spaces closures are being taken.
The intersection of these sets over all $t_0$ yields
\[
  \bigcap_{t_0}\cl\Big(\bigcup_{t\ge t_0}\Gamma_{t\cdot s}\Big)\ =\ 
  P_l \Big(\overline{\Gamma_{s_l}}\cap (F\times\{0\})\Big)\ =\ C(s)\,.
\]
For the second point let $K,U\subset F$ be as stated in the hypothesis. Consider
the other projection $Q_l: F\times \rz^l\rightarrow \rz^l$. The restriction of $Q_l$ to
$\overline{\Gamma_{s_l}}\cap P_l^{-1}(K)$ is proper, and for the fiber over
$0\in \rz^l$ the following inclusion holds:
\[
  \overline{\Gamma_{s_l}} \cap(K\times\{0\})\subset P_l^{-1}(U)\,.
\]
Then the same inclusion holds for $\overline{\Gamma_{s_l}}\cap
(K\times\{v\})$, $v\in B_{t_0^{-1}}(0)$, for some $t_0>0$. But this means
\[
  \Gamma_{t\cdot s}\cap K\subset U\quad \forall\,t>t_0\,.
\]
2)\ One inclusion follows from $\Gamma_{s_l}\supset Z(s)\times
(\rz^l\setminus \{0\})$. Conversely, let $f\in F_x$ and $s(x)\neq 0$. Let $K$ be a
compact neighbourhood of $f$ in $F$ such that $s$ has no zeros on $q(K)$.
Then
\[
  A\ =\ \{(t,y)\in \rz_{\ge0}\times N\mid t\cdot s(y)\in K\}
\]
is compact. Choosing $c>0$ with $\prj_1(A) \subset [0,c]$ we obtain
\[
  \Gamma_{t\cdot s}\cap K\ =\ \emptyset\quad\mbox{for any $t>c$}\,.
\]
Hence $f\not\in C(s)$.
\qed

The reason for introducing the factor $\rz^l$ with $l>1$ is the long
exact sequence of homology groups (of the second kind, that is with
locally finite singular chains) 
\[
  H_{n+l}(C(s))\longrightarrow H_{n+l}(\overline{\Gamma_{s_l}})
  \longrightarrow H_{n+l}(\Gamma_{s_l})\longrightarrow H_{n+l-1}(C(s))\,.
\]
Provided $n+l-1>\dim C(s)$ (e.g.\ $l>\dim F-n+1$) the groups on the left and
right vanish by the general vanishing theorem of homology, cf.\ e.g.\
\cite[IX.1.Prop.1.6]{iversen}. So the (orbifold!) fundamental class
\[
  [\Gamma_{s_l}]\ =\ (s_l)_*[N\times\rz^l]
\]
extends uniquely to an $(n+l)$-dimensional homology class on
$\overline{\Gamma_{s_l}}$, that will be conveniently denoted
$[\overline{\Gamma_{s_l}}]$ (slight abuse of notation). Here we assumed
(without loss of generality) $N$ to be pure $n$-dimensional and we chose an
orientation of $\rz^l$ that will finally drop out. We are now in position to
define the limit of $[\Gamma_{t_s}]$ as $t$ tends to $\infty$.
\begin{defi}\rm
Let $\delta_0\in H_{\{0\}}^l(\rz^l)$ be Poincar\' e-dual to $\{0\}\subset \rz^l$ with
respect to the chosen orientation. Let $Q_l: p_l^*F\rightarrow \rz^l$ be the
projection. Noticing that $C(s)=\overline{\Gamma_{s_l}}\cap Q_l^{-1}(0)$ we
define
\[
  [C(s)]\ :=\ [\overline{\Gamma_{s_l}}]\cap Q_l^*\delta_0 \in H_n(C(s))\,.
\]
\vspace{-8ex}

\qed
\end{defi}
Implicit in the notation is the first statement in
\begin{prop}
  $[C(s)]$ is independent of $l$ and homologous to $[\Gamma_{t\cdot s}]$ for
  any $t\neq 0$ as class on $F$.
\end{prop}
\pf
We identify $\rz^l$ with a linear subspace in $\rz^{l+1}$. Let $\eta\in H^1_{\rz^l}
(\rz^{l+1})$ be the corresponding cohomology class with supports.
We write $\delta^k$ for the integral generator of $H^k_{\{0\}}(\rz^k)$ (previously
denoted $\delta_0$). Then
\[
  [\Gamma_{s_l}]\ =\ [\Gamma_{s_{l+1}}]\cap Q_{l+1}^*\eta,\quad
  \delta^{l+1}\ =\ \eta\cup\delta^l\,,
\]
hence
\begin{eqnarray*}
  [\overline{\Gamma_{s_l}}]\cap Q_l^*\delta^l&=&
  \Big([\overline{\Gamma_{s_{l+1}}}]\cap Q_{l+1}^*\eta\Big)
  \cap Q_l^*\delta^l\\
  &=&[\overline{\Gamma_{s_{l+1}}}]\cap
  Q_{l+1}^*(\eta\cup\delta^l)\ =\ 
  [\overline{\Gamma_{s_{l+1}}}]\cap Q_{l+1}^*\delta^{l+1}\,.
\end{eqnarray*}
And for $t>0$ the Poincar\' e-dual $\delta_t$ to $\{t\}\subset\rz^l$ is
cohomologous to $\delta^l$ as class on $\rz^l$. Thus
\[
  [\Gamma_{t\cdot s}]\ =\ [\Gamma_{s_l}]\cap Q_l^*\delta_t\ =\ 
  [\overline{\Gamma_{s_l}}]\cap Q_l^*\delta_t\ =\
  [\overline{\Gamma_{s_l}}]\cap Q_l^*\delta^l\ =\ [C(s)]
\]
in $H_n(F)$.
\qed
\vspace{2ex}

In a holomorphic situation we retrieve the following familiar picture
\cite[\S14.1]{fulton}: Let $F$ be a holomorphic vector bundle over a
complex manifold $N$, and let $Z$ be the zero locus of a holomorphic section
$s$ of $F$. The differential of $s$ induces a closed embedding $\iota$ of the
{\em normal bundle} $N_{Z|N}$ of $Z$ in $N$ into $F$. $N_{Z|N}$ is the linear
fiber space over $Z$ associated to the conormal sheaf $\cali/\cali^2$,
$\cali$ the ideal sheaf of $Z$ in $N$. The {\em normal
cone} $C_{Z|N}$ (the analytic analog of $\Spec_Z \oplus_{d\ge0} \cali^d
/\cali^{d+1}$) is a closed subspace of $N_{Z|N}$. Thus $\iota(C_{Z|N})$ is a closed
subcone of $F$.
\begin{prop}\label{analyt_cone}
  $C(s)=\iota(C_{Z|N})$ and $[C(s)]=\iota_*[C_{Z|N}]= [\iota(C_{Z|N})]$
  (where $[C_{Z|N}]$, $[\iota(C_{Z|N})]$ are the fundamental classes
  of the correponding complex spaces).
\end{prop}
\pf
Our construction of $C(s)$ and $[C(s)]$ is nothing but (a real version of) the
``deformation to the normal cone'', which in this case states that
$\iota(C_{Z|N})$ can be obtained as analytic limit of
$\Gamma_{t\cdot s}$, $t\in\cz$, $|t| \rightarrow \infty$ \cite[Rem.5.1.1]{fulton}.
From this the compatibility of the two limits can be easily deduced.
\qed

\noindent
It should be clear that the same conclusions hold in the category of analytic
orbifolds, but we will not need this.
\vspace{2ex}

For later reference we also observe here two simple lemmata:
\begin{lemma}\label{unif_cones}
  Let $q: \hat F_U=\hat U\times\rz^r\stackrel{/G^F}{\rightarrow} F|_U$ be a
  local uniformizing trivialization of a topological orbibundle $F$ over a local
  uniformizing system $U=\hat U/G$ of the base orbifold $N$, and let $s$ be a
  section of $F$ uniformized by $\hat s_U: \hat U\rightarrow \hat F_U$. Let $b$
  be the generic covering degree of $q$. Then
  \[
    C(s)\ =\ q(C(\hat s))\,,\quad
    [C(s)]\ =\ \frac{1}{b} q_*[C(\hat s)]\,.
  \]
\end{lemma}
\pf
This follows immediately by the corresponding identities for $\Gamma_{t\cdot
s}$ and $\Gamma_{t\cdot \hat s}$.
\qed


\subsection{Local decomposition of cone classes}
Returning to GW-theory we will choose in Section~4.1 a morphism $\tau:
F\rightarrow \cale$ such that locally there is complex subbundle $F^h\subset
F$ with $\tau^h=\tau|_{F^h}$ (essentially) holomorphically spanning the
cokernel of the linearization of $s=s_\dbar$, cf.\ Proposition~\ref{tau}.
We obtain two finite-dimensional oriented orbifolds $\tilde Z\subset F$ and
$\tilde Z^h\subset F^h$ as zero loci of $\tilde s=q^*s+\tau$ and $\tilde
s^h=(q^h)^*s +\tau^h$ respectively, $q:F\rightarrow \calc(M,p)$,
$q^h:F^h\rightarrow \calc(M,p)$ the projections. The tautological sections
$s_\can$ and $s_\can^h$ of $q^*F$ and $(q^h)^*F$ both have zero locus
$\calc(M)$. The associated cones and cone classes defined in the last section
will be written
\begin{eqnarray*}
  C(\tau)&=& C(s_\can)\ \subset F|_Z,\ [C(\tau)]\\
  C(\tau^h)&=& C(s_\can^h)\ \subset F^h|_Z,\ [C(\tau^h)]\,.
\end{eqnarray*}
For a decomposition $F=F^h\oplus\bar F$ and homology classes $\alpha,\beta$
supported on closed subsets $A\subset F^h$, $B\subset\bar F$ let us write
\[
  \alpha\oplus\beta\ :=\ (\alpha\times\beta)\cap
  (q^h\times\bar q)^*\delta_\Delta
\]
for their {\em direct sum}, where $q^h\times\bar q: F^h\times\bar F
\rightarrow N\times N$ is the product of the bundle projections and
$\delta_\Delta \in H_\Delta^*(N\times N)$ is Poincar\' e-dual to the diagonal
$\Delta\subset N\times N$.

The object of this section is to show that $C(\tau)$ is already determined by
$C(\tau^h)$.
\begin{prop}\label{cone_decomp}
  Let $F=F^h\oplus\bar F$ be a decomposition into complex orbibundles in such
  a way that
  \begin{itemize}
  \item
    $\tau^h:=\tau|_{F^h}$ is injective and spans the cokernel of the linearization
    $\sigma$ relative a local, finite dimensional parameter space along
    $\calc(M)=Z(s)$ and has the regularity properties of $\tau$
    (cf.\ \cite[Def.1.15]{si1})
  \item
    $\bar\tau:= \tau|_{\bar F}$ maps to $\img\sigma$ along $\calc(M)$.
  \end{itemize}
  Then $C(\tau)=C(\tau^h) \oplus\bar F$ and $[C(\tau)]=[C(\tau^h)] \oplus
  [\bar F]$.
\end{prop}
Before turning to the proof three remarks are in order: First, while $s=s_\dbar$ is
not in general globally differentiable, locally it is so relative to a
parameter space $S$ of a semiuniversal deformation of the curve.
The corresponding relative differential is nothing but the linear
$\dbar$-operator from ${\check L}_1^p(C,\ph^*T_M)$ to ${\check L}^p(C,
\ph^*T_M \otimes\bar\Omega)$ (restricted to a complement
of the kernel in case $(C,{\bf x})$ is not stable). Second, it will be crucial that
differentiability properties are imposed only on $\tau^h$. $\bar F$ will indeed
only be constructed as topological subbundle. And third, the proposition
together with Proposition~\ref{analyt_cone}
shows that, locally, $C(\tau)$ is the product of an analytic cone over $\calc(M)$,
pure-dimensional of dimension equal to expected dimension
(index of relative differential plus $\dim S$ minus $\dim \Aut(C,{\bf x})$)
plus $\rk F^h$, and a complex vector space of dimension $\rk\bar F$.
\vspace{1ex}

\noindent
\pf
Let $\bar s_\can$ be the tautological section of $q^*\bar F$. Consider the
section $\theta$ of $q^*F\times\rz^l\times\rz \rightarrow \tilde
Z\times\rz^l\times\rz$ over $F \times
(\rz^l\setminus\{0\})\times (\rz\setminus\{0\})$ defined by
\[
  \theta(z,v,u)\ :=\ |v|^{-1}\cdot (s_\can^h+u\cdot\bar s_\can)\,.
\]
This should be viewed as a two-parameter family of sections of $q^*F$
interpolating between $(s_\can)_l= |v|^{-1}\cdot s_\can$ ($u=1$) and
$(s_\can^h)_l\oplus 0$ ($u=0$). Write $\theta_u$ for the restriction of $\theta$
to $q^* F\times\rz^l \times\{u\}$ and $Q_F$, $Q_l$, $Q$ for the projections
from $q^*F \times \rz^l\times\rz$ to the three factors. We denote by
$\delta_0^{\rz^l}\in H_{\{0\}}^l (\rz^l)$ and $\delta_u^\rz\in
H^1_{\{u\}}(\rz)$ the Poincar\'e duals to $\{0\}\subset \rz^l$, $\{u\}\subset \rz$.
By definition
\[
  [C(\tau)]\ =\ [\overline{\Gamma_{\theta_1}}]\cap Q_l^*\delta_0^{\rz^l}\,,
\]
as classes supported on $C(\tau) \times\{0\}\times\{1\} \subset
q^*F\times\{0\} \times \{1\}$. $[\overline{ \Gamma_{\theta_1}}] $ is the unique
class extending the fundamental class of the oriented orbifold
$\Gamma_{\theta_1}$. In the construction of local uniformizers for
$\Gamma_\theta$ we may take $u$ as parameter in the application of the
implicit function theorem, cf.\ the proof of \cite[Thm.1.21]{si1}. This shows
\[
  [\Gamma_{\theta_1}]\ =\ [\Gamma_\theta]\cap Q^*\delta_1^\rz
\]
as classes supported on $\Gamma_{\theta_1}$. Since by Lemma~\ref
{limit_product} below $\overline{\Gamma_{\theta_u}}= \overline{
\Gamma_\theta} \cap Q^{-1} (u)$ (the disjoint union of $\Gamma_{
\theta_u}$ and $C(\tau)\times \{0\}\times \{u\}$), we get at $u=1$
\[
  A_1\ :=\ [C(\tau)]\times \{0\}\times \{1\}\ =\ [\overline{\Gamma_\theta}]
  \cap (Q_l\times Q)^*(\delta_0^{\rz^l}\times \delta_1^\rz)\,.
\]
And by the same lemma
\[
  A_u\ :=\ [\overline{\Gamma_\theta}]
  \cap (Q_l\times Q)^*(\delta_0^{\rz^l}\times \delta_u^\rz),\quad u\in\rz
\]
is a family of classes on $C(\tau)\times \{0\}\times \rz$, and these are
mutually homologic because the $\delta_u^\rz$ are cohomologous classes
on $\rz$.

Together with $\overline{\Gamma_{\theta_0}}\cap Q_l^{-1}(0)= (C(\tau)\oplus
\bar F)\times \{0\}\times \{0\}$ we obtain
\[
  [C(\tau)]\ =\ (Q_F)_*A_1\ =\ (Q_F)_*A_0\ =\ 
  [C(\tau^h)]\oplus [\bar F]\,.
\]
The set-theoretic part of the claim is also proved in the lemma below.
\qed

We still owe the set-theoretic part of the lemma:
\begin{lemma}\label{limit_product}
  1)\ \ $\overline{\Gamma_\theta}\cap (q^*F\times \{0\}\times\rz)
  \ =\ C(\tau)\times\{0\}\times \rz$\\[1ex]
  2)\ \ $C(\tau)=C(\tau^h)\oplus \bar F$ (set-theoretically).  
\end{lemma}
\pf
We will show the inclusions
\[
  \overline{\Gamma_\theta}\cap Q_l^{-1}(0)\subset
  (C(\tau^h)\oplus \bar F)\times \{0\}\times\rz,\quad
  C(\tau^h)\oplus\bar F\subset C(\tau)\,.
\]
The lemma will then be finished with $C(\tau)=\overline{\Gamma_{\theta_1}}
\cap Q_l^{-1}(0)$ (hence (2)) and $C(\tau^h)\oplus \bar F=
\overline{\Gamma_{\theta_0}} \cap Q_l^{-1}(0)$, together with the
observation that rescaling $(f^h,\bar f,v,u) \mapsto (f^h,u\cdot \bar f, v,1)$
gives
\[
  \overline{\Gamma_\theta}\cap\Big(q^*F\times\{0\}\times (\rz\setminus\{0\})\Big)
  \ =\ (C(\tau^h)\oplus\bar F)\times \{0\}\times(\rz\setminus\{0\})\,.
\]

In view of Proposition~\ref{cones_prop} we have to prove the following: Let
$t_\nu\in \rz$, $u_\nu\in\rz$, $f_\nu\in F_{x_\nu}$ be sequences with
\begin{eqnarray}\label{def_eqn}
  s_\dbar(x_\nu)-t_\nu\cdot\Big(\tau_{x_\nu}^h(f_\nu^h)
  -u_\nu\cdot \bar\tau_{x_\nu} (\bar f_\nu)\Big)\ =\ 0
\end{eqnarray}
(cf.\ the definition of $\theta$ together with that of $\tilde Z \subset F$ as zero
locus of $q^* s_\dbar-\tau$) and
\[
  t_\nu\longrightarrow 0,\quad u_\nu\longrightarrow u,\quad
  f_\nu\longrightarrow f\,.
\]
As before we use superscript ``h'' and a bar to denote components in $F^h$ and
$\bar F$ respectively and $\tau_y$ (etc.) to denote the restriction of $\tau$ to
$F_y$ . Any $(f,0,u)\in \overline{\Gamma_\theta}\cap
Q_l^{-1}(0)$ is of this form. We claim $f^h\in C(\tau^h)$.

To this end we want to work on adapted charts. Recall that at $z=(C,{\bf x},
\ph)$, $s_\dbar$ was locally uniformized by
\[
  \hat s_\dbar:\ S\times V\longrightarrow E_0\,,
\]
$S$ the parameter space of deformations of the domain of the curve, $V\subset
{\check L}_1^p (C;\ph^*T_M)$ (an open set in a linear subspace) of finite
codimension and $E_0= {\check L}^p(C; \ph^*T_M\otimes \bar\Omega)$
uniformizing $\cale_z$. This map was differentiable relative
$S$ (i.e.\ for fixed $s\in S$) with relative differential $D_V \hat s_\dbar$
uniformly continuous at the center $(0,0)\in S\times V$. Write $\sigma
= D_V \hat s_\dbar (0,0)$. By the regularity
properties of $\hat\tau^h$ and since $\hat\tau_z^h$ is injective and spans
$\coker\sigma$ we may change the trivialization of $\cale$ in such a way that
${\hat F}^h$ is identified via ${\hat \tau}^h$ with its image $C$ on the
central fiber $E_0$.

Choose a complementary subspace $P\subset V$ to $K:=\ker \sigma$, set
$Q=\sigma(P)$, and write $\prj_Q:E_0 \rightarrow Q$, $\prj_C: E_0\rightarrow
C$ for the projections with kernel $C$ respectively $Q$. We can now apply the
implicit function theorem to $\prj_Q\circ \hat s_\dbar: S\times K\times P
\rightarrow Q$ with parameter space $W=S\times K$. We can thus change
coordinates on $S\times V$ in such a way that
\begin{eqnarray*}
  \hat s_\dbar (w,p)&=&(p,\kappa(w,p))\ \in\ Q\times C\\
  {\hat \tau}_{(w,p)}^h(f^h)&=&(0,f^h)\\
  {\hat {\bar\tau}}_{(w,p)}(\bar f)&=& (\prj_Q {\hat{\bar\tau}}_{(w,p)}(\bar f),
  \prj_C{\hat{\bar\tau}}_{(w,p)}(\bar f))
\end{eqnarray*}
with $\kappa: W\times P\rightarrow C$ differentiable relative $W$,
$D_P\kappa(0,0) =0$, $D_P\kappa$ uniformly continuous at $(0,0)$, and
$\prj_C\bar\tau_{(0,0)} =0$ ($\img {\hat{\bar\tau}}_{(0,0)} \subset
\img\sigma$ by hypothesis).

Since the structure map $S\times V\rightarrow \calc(M;p)$ is locally proper, it
suffices to prove the claim about the limit on local uniformizers. For readabilties
sake we will drop the hats that usually indicate local unifromizers. Write
$x_\nu=(w_\nu, p_\nu)\in W\times P$, $x=(0,0)$. Equation~\ref{def_eqn} now
splits into the two equations (in $Q$ and $C$ respectively)
\begin{eqnarray*}
  p_\nu-t_\nu u_\nu \prj_Q\bar\tau_{(w_\nu,p_\nu)}(\bar f_\nu)&=&0\\
  \kappa(w_\nu,p_\nu)-t_\nu f_\nu^h-t_\nu u_\nu\prj_C
  \bar\tau_{(w_\nu,p_\nu)}(\bar f_\nu)&=&0\,.
\end{eqnarray*}
We claim that
\[
  f'_\nu\ :=\ \Big(\frac{1}{t_\nu}\kappa(w_\nu,0),\bar f_\nu\Big)
\]
has the same limit $(f^h,\bar f)$ as $f_\nu$. From the first equation we see
that ($t_\nu\neq 0$)
\[
  \frac{p_\nu}{t_\nu}\longrightarrow u\cdot \prj_Q\bar\tau_{(0,0)}(\bar f)\ 
  =\ u\cdot \bar\tau_{(0,0)}(\bar f)\,.
\]
By uniform continuity of $D_P\kappa$ and because $D_P\kappa(0,0) =0$
\[
  \frac{1}{|t_\nu|}\Big|\kappa(w_\nu,p_\nu)-\kappa(w_\nu,0)\Big|\ \le\ 
  |\!|D_P\kappa|\!|_{B_{|p_\nu|}(0)}\cdot\Big|\frac{p_\nu}{t_\nu}\Big|\ 
  \longrightarrow\ 0\,.
\]
Together with
\[
  \frac{1}{t_\nu}\kappa(w_\nu,p_\nu)\ =\ 
  f_\nu^h+ u_\nu\cdot\prj_C\bar\tau_{(w_\nu,p_\nu)}(\bar f_\nu)\ 
  \longrightarrow\ f^h
\]
(which is where the assumption $\img\bar\tau\subset \img\sigma$ comes in)
this establishes the claim and hence the first inclusion, for $(f'_\nu,t_\nu)$ are
in $\Gamma_{\theta_0}$.

Turning to the second inclusion $C(\tau^h)\oplus\bar F\subset C(\tau)$ we
replace $\bar F$ by a subbundle ${\bar F}'\subset F$ with ${\bar F}'_0 =\bar
F_0$ and such that ${\bar\tau}'
=\tau|_{{\bar F}'}$ has the regularity properties of $\tau$. To $f^h\in C(\tau^h)$
choose sequences $0\neq t_\nu \rightarrow 0$, $F_{x_\nu= (w_\nu,0)} \ni
f_\nu^h \rightarrow f^h$ with
\[
  \kappa(w_\nu,0)-t_\nu f_\nu^h\ =\ 0\,.
\]
For any $\bar f\in \bar F_0 =\bar F'_0$ we want to find ${f_\nu^h}'\rightarrow f^h$
with
\begin{eqnarray*}
  p_\nu-t_\nu \prj_Q{\bar\tau}'_{(w_\nu,p_\nu)}(\bar f)&=&0\\
  \kappa(w_\nu,p_\nu)-t_\nu {f_\nu^h}'-t_\nu \prj_C
  {\bar\tau}'_{(w_\nu,p_\nu)}(\bar f)&=&0\,.
\end{eqnarray*}
To the first equation we may apply the implicit function theorem with
parameters $t_\nu,w_\nu$ to conclude unique existence of $p_\nu
\rightarrow 0$ for $\nu$ large from the solution $p=0$, $t=0$, $w=0$ in the
limit $\nu\rightarrow \infty$. The second equation in turn forces
\[
  {f_\nu^h}'\ =\ \frac{1}{t_\nu}\kappa(w_\nu,p_\nu)-
  \prj_C{\bar\tau}'_{(w_\nu,p_\nu)}(\bar f)\,.
\]
Since as above
\[
  {f_\nu^h}'-f_\nu^h\ =\ \frac{1}{t_\nu}\Big(\kappa(w_\nu,p_\nu)-
  \kappa(w_\nu,0)\Big)\ \longrightarrow\ 0
\]
we deduce ${f_\nu^h}'\rightarrow f^h$ as claimed.
\qed


\section{Comparison of algebraic and limit cone}
\subsection{Choice of Kuranishi structure}
Let $(M,\omega)$ be a K\"ahler manifold. Recall the construction of \cite{si1}
applied to $(M,\omega)$ viewed as symplectic manifold with almost complex
structure the integrable one: The space $\calc (M;p) =\bigcup_ {R,g,k} \calc_
{R,g,k} (M;p)$ of stable marked complex curves in $M$ of Sobolev class $L_1^p$
is a Banach orbifold. The map $(C,{\bf x}, \ph) \mapsto \dbar \ph$ is a section
$s_\dbar$ of the Banach orbibundle $\cale$ with fibers $\cale_{(C,{\bf x},\ph)}=
{\check L}^p (C;\ph^* T_M\otimes \bar\Omega_C)$. Local uniformizing systems
at $(C,{\bf x},\ph)$ are of the form $S\times V$ with $S$ the base of an
analytically semiuniversal deformation of $(C,{\bf x})$ and $V$ an open subset
of a linear subspace of ${\check L}_1^p (C;\ph^*T_M)$ of codimension equal to
$\dim\Aut (C,{\bf x})$. In such a chart $\hat s_\dbar$ is differentiable relative $S$,
with differential $\sigma$ relative $S$ a family of Fredholm operators that is
uniformly continuous at $(0,0) \in S\times V$. By spanning the cokernel of
$\sigma$ along some compact part of $\calc(M)$, say $\calc_{R,g,k} (M)
=\calc(M)\cap \calc_{R,g,k} (M;p)$, by sections of $\ph^*T_M\otimes
\bar\Omega_C$ supported away from the singularities of $C$, parallel transport
by local trivializations of $\cale$ and multiplication with a bump function in a
neighbourhood of $\calc_{R,g,k} (M;p)$ in $\calc(M;p)$ we constructed a
morphism (``Kuranishi structure'')
\[
  \tau: F\longrightarrow \cale
\]
from a finite rank complex orbibundle $F$ living on a neighbourhood of
$\calc_{R,g,k} (M;p)$ in $\calc(M;p)$ with the following properties:
\begin{itemize}
\item
  $\tau$ has the same differentiability properties as $s_\dbar$
\item
  $\tau$ spans the cokernel of the linearization of $s_\dbar$ along $Z(s_\dbar)$,
  i.e.\ for any $(C,{\bf x}, \ph)\in\calc(M)$
  \[
    \img \hat\tau_{(C,{\bf x},\ph)} + \img\sigma_{(C,{\bf x},\ph)}\ =\ 
    \hat\cale_{(C,{\bf x},\ph)}\,.
  \]
\end{itemize}
The section $\tilde s:= q^*s-\tau$, $q:F\rightarrow \calc(M;p)$ the bundle
projection, is then a transverse (locally relative $S$) section of the
Banach orbibundle $q^*\cale$ over the total space of $F$.
$\tilde Z=Z(\tilde s) \subset F$ is thus an oriented, finite-dimensional orbifold.
Let $\Theta_F \in H^{\rk F}_{\calc(M;p)} (F)$ be the Thom class of $F$
($\calc(M;p)$ identified with the zero section of $F$). The virtual fundamental
class of $\calc_{R,g,k}(M)$ was defined by
\[
 \GW_{R,g,k}^M\ :=\ [\tilde Z]\cap \Theta_F\,.
\]
Note that the restriction of $\Theta_F$ to $\tilde Z$ can also be written
$s_\can^*\Theta_{q^* F}$, so this definition is essentially finite dimensional.
\vspace{2ex}

To compare with the complex analytic definition of virtual fundamental classes,
as a first try one might want to make the whole symplectic construction Banach
analytic. This seems to be hard if not impossible. It is however easy
to gain enough analyticity locally to make the comparison with the analytic
treatment given in the first chapter work.

At holomorphic $(C,{\bf x},\ph)$ let $\Theta:\calb\rightarrow \prod_{j\le m'}
L^p(\Delta^j;\cz^n)$ be a holomorphic Kuranishi model as constructed in
Section~2.3. Let $\breve\pi: \breve \Gamma:= \calb\times _S \calc$
be the universal curve over $\calb$ and $\breve \Phi:\breve\Gamma
\rightarrow M$ be the evaluation map (generally we will use
the accent \ $\breve{}$\  for objects on $\calb$).
We will write $U_0\subset\breve \Gamma$ for the union of the open sets
formerly denoted $\Sigma$, $U_i$ and $W_j$ for $j>m'$, and $U_i$ for $W_i$,
$i=1,\ldots, m'$. So $\{U_i\}$ is an open covering of $\breve\Gamma$ that is
Stein relative $\calb$. Let $\breve\calt$ be the tangent bundle of $\calb$
relative $S$ and $\breve\cale:=\calb\times\prod_{j\le m'} L^p(\Delta^j;\cz^n)$.
The former has fiber $\check L_1^p(C';{\ph'}^* T_M)\cap \calo^{(c)}
(U_0;\ph^*\calt_M)$ at $(C',{\bf x'},\ph')\in\calb$ while the latter should be viewed
as a version of the Banach bundle $\cale\downarrow\calc(M;p)$ on $\calb$.
Note that while $\breve\calt$ and $\breve\cale$ parametrize non-holomorphic
objects they are {\em holomorphic} Banach bundles over $\calb$. The Fredholm
map $\Theta$ exhibiting $\Hom_S(\calc,M)$ as fiber over $0$ can now be viewed
as holomorphic section of $\breve\cale$.

Similarly, while the evaluation map is not holomorphic (along the
fibers of $\breve\pi$) $\breve\Phi^* T_M$ is a holomorphic vector bundle
over $\breve\Gamma$, local holomorphic trivializations being given
by pull-back of a frame of local holomorphic vector fields on $M$.

While $U_0$ might now have singularities a straightforward modification
of the arguments in \cite[\S4.2]{si1} shows that the spaces of
relative holomorphic \v Cech cochains $\breve\pi_*^{(i)}\breve\Phi^*
\calt_M$ are still holomorphic Banach bundles over $\calb$ (and the
same holds true by replacing $\breve \Phi^* \calt$ by any finite rank
holomorphic vector bundle over $\breve\Gamma$). Recall also from op.cit.\
that (the restrictions to $\calb\subset\calc(M;p)$) of the tangent bundle
of $\calc(M;p)$ relative $S$ and the Banach bundle $\cale$ can be written
\[
  \calt\ =\ \mbox{${{\breve\pi}_1^p}$}_*\breve\Phi^*\calt_M\,,\quad
  \cale\ =\ \breve\pi^p_*(\breve\Phi^*\calt_M\otimes_\cz
  \bar\Omega_{\breve\Gamma/\calb})\,.
\]
We will also need to extend the holomorphic bundles $G$ and $H$ from
Section~1.2 to $\calb$. The sheaf of sections of $H$ was $R^1\pi_* \caln$,
where $\caln= \pi^* \pi_*(\Phi^* \calt_M\otimes \call^\nu) \otimes
\call^{\otimes -\nu}$, $\call= \omega (x_1+\ldots +x_k)\otimes
\Phi^*\calh_M^{\otimes 3}$, fits into an exact sequence of locally
free sheaves
\[
  0\longrightarrow\calk \stackrel{\kappa}{\longrightarrow} \caln
  \stackrel{\nu}{\longrightarrow} \Phi^*\calt_M\longrightarrow 0\,.
\]
Now $\Phi^*T_M$ and $\call$ extend naturally holomorphically to
$\breve\Gamma$ and so do $\caln$ and $\calk$ and the above sequence.
Let $\breve\caln$ and $\breve \calk$ denote these extensions.
We thus obtain a commutative diagram of holomorphic Banach bundles
over $\calb$ with exact rows and columns
\begin{eqnarray}\label{big_diag}&\begin{array}{ccccccccc}
&&0&&0\\[3pt]
&&\longdownarrow&&\longdownarrow\\[3pt]
0&\llongrightarrow&\breve\pi_*^{(0)}\breve\calk
&\stackrel{\kappa^{(0)}}{\llongrightarrow}& \breve\pi_*^{(0)}\breve\caln
&\stackrel{\nu^{(0)}}{\llongrightarrow}& \breve\pi_*^{(0)} \breve\Phi^*\calt_M
&\llongrightarrow &0\\[3pt]
&&\longdownarrow&&\longdownarrow&&\longdownarrow\\[3pt]
0&\llongrightarrow&\breve\pi_*^{(1)}\breve\calk
&\stackrel{\kappa^{(1)}}{\llongrightarrow}& \breve\pi_*^{(1)}\breve\caln
&\stackrel{\nu^{(1)}}{\llongrightarrow}& \breve\pi_*^{(1)} \breve\Phi^*\calt_M
&\llongrightarrow &0\\[3pt]
&&\diagl{q_{\breve G}}&&\diagr{q_{\breve H}}\\[3pt]
&&\breve G&\stackrel{R^1\breve\pi_*\kappa}{\llongrightarrow}&
\breve H\\[3pt] &&\longdownarrow&&\longdownarrow\\[3pt]
&&0&&0
\end{array}\end{eqnarray}
where we {\em define} $\breve G$ and $\breve H$ as cokernels of the first two
columns. These are holomorphic extensions of $G$ and $H$ to a neighbourhood
of $(C,{\bf x},\ph)$ in $\calb$.

Recall that we called a commutative square of Banach bundles
\[\begin{array}{ccc}
  E&\stackrel{\gamma}{\llongrightarrow} &G\\[3pt]
  \diagl{\alpha}&&\diagr{\beta}\\[3pt]
  F&\stackrel{\delta}{\llongrightarrow}& H
\end{array}\]
a {\em quasi-isomorphism} (between $\alpha$ and $\beta$, and between
$\gamma$ and $\delta$) iff the sequence
\[
  0\longrightarrow E\stackrel{(\alpha,\gamma)}{\longrightarrow}
  F\oplus G\stackrel{\delta-\beta}{\longrightarrow} H\longrightarrow 0
\]
is exact \cite[Def.4.8]{si1}. Equivalently, $(\alpha,\beta)$ induces
fiberwise isomorphisms between kernels and cokernels
of $\gamma$ and $\delta$ (or the other way
around). In op.cit.\ we also required this sequence to be split in the case
of which the square is not only cartesian ($E\simeq F\oplus_H G$) but
also cocartesian ($H\simeq (F\oplus G)/ E$, as Banach bundles!). All our
quasi-isomorphisms will in fact be split but since we will never need this
property we will not verify it.

Now a diagram like (\ref{big_diag}) above always induces a
quasi-isomorphism between $[\breve G\rightarrow \breve H]$ and
$[\breve\pi_*^{(0)}\breve \Phi^*\calt_M \rightarrow \breve\pi_*^{(1)}
\breve \Phi^*\calt_M]$
(unique up to homotopy) locally as follows: Let $\eta$ and $\theta$
be (local) holomorphic right-inverses to $q_{\breve G}$ and $q_{\breve H}$.
Then
\[
  q_{\breve H}(\theta\circ R^1\breve\pi_*\kappa-\kappa^{(1)}\circ \eta)
  \ =\ 0\,,
\]
so the term in the bracket lifts (uniquely) to a map
\[
  \zeta:\ \breve G\ \longrightarrow\ \breve\pi_*^{(0)}\caln\,.
\]
Define
\[
  \alpha^{(0)}\ :=\ \nu^{(0)}\circ\zeta,\quad
  \alpha^{(1)}\ :=\ \nu^{(1)}\circ\theta\,.
\]
Recall also the morphisms \cite[\S4.3]{si1}
\[\begin{array}{rclrcl}
  \breve\pi_*^{(0)}\breve\Phi^*\calt_M
  &\longrightarrow& \calt\,,\quad&
  (v_i)&\longmapsto&\sum_i \rho_i v_i\\[3pt]
  \breve\pi_*^{(1)}\breve\Phi^*\calt_M
  &\longrightarrow& \cale\,,\quad&
  (v_{ij})&\longmapsto&\frac{1}{2}\sum_{i,j} v_{ij}\cdot\dbar\rho_i\,.
\end{array}\]
Here $\rho_i$ is a partition of unity subordinate to $U_i$, so these maps
indeed factor over the inclusions $\breve\calt\hookrightarrow \calt$
and $\breve\cale\hookrightarrow\cale$.
\begin{lemma}\label{GH_quasiiso}
  The squares in the following diagram are quasi-isomorphisms.
  \[\begin{array}{ccccccc}
    \breve G&\stackrel{\alpha^{(0)}}{\longrightarrow}&
    \breve\pi_*^{(0)}\breve\Phi^*\calt_M
    &\stackrel{\beta^{(0)}}{\longrightarrow}&
    \breve\calt&\longrightarrow&\calt\\[2pt]
    \longdownarrow&&\longdownarrow&&
    \longdownarrow&&\diagr{\sigma=\dbar}\\[3pt]
    \breve H&\stackrel{\alpha^{(1)}}{\longrightarrow}&
    \breve\pi_*^{(1)}\breve\Phi^*\calt_M
    &\stackrel{\beta^{(1)}}{\longrightarrow}&
    \breve\cale&\longrightarrow&\cale
  \end{array}\]
\end{lemma}
\pf
The claim for the first square follows from chasing Diagram~\ref{big_diag}.
That the last square induces isomorphisms on kernels and cokernels of the
vertical maps has been proven in Proposition~\ref{Tbar_quasiiso}.
The composition of the right two squares is a quasi-isomorphism
essentially by \cite[Cor.4.11]{si1}, observing that there is a (family of)
right-inverses $T$ to the $\dbar$-operator on $U_0$ (cf.\ Section~2.1).
This shows that also the middle square is a quasi-isomorphism.
\qed

Notice that $F$ extends naturally as finite rank
orbibundle to a neighbourhood of $\calc_{R,g,k}(M)$ in
$\calc(M;p)$. In fact, the finite rank orbibundle in
\cite[\S6.4]{si1} was a direct sum of bundles of this form. Locally we may
thus take for $\tau$ an extension of $\beta^{(1)}\circ \alpha^{(1)}$ to a
local uniformizing system of $\calc(M;p)$, multiplied by a bump function,
as specified in \cite[\S6.5]{si1}.
\begin{prop}\label{tau}
  For any $R\in H_2(M;\gz)$, $g,k\ge 0$, the Kuranishi structure $\tau:
  F\rightarrow \cale$ for $s_\dbar$ may be chosen in such a way that
  \begin{itemize}
  \item
    $F=\oplus_{\nu=1}^l F_\nu$ with each $F_\nu$ restricting to the orbibundle
    $H$ along $Z$
  \item
    for any $(C,{\bf x},\ph) \in\calc_{R,g,k}(M)$ there exists an open
    neighbourhood $U\subset \calc(M;p)$ and a $\nu$ with $\tau_\nu=
    \tau|_{F_\nu}$ an extension of $\beta^{(1)}\circ\alpha^{(1)}$. In
    particular,
    \[
      \tilde Z_\nu\ :=\ Z(q_\nu^*s_\dbar+\tau_\nu)\ \subset\ F_\nu
    \]
    is a complex orbifold at $(C,{\bf x},\ph)$.
  \end{itemize}
\end{prop}
\pf
$\tilde Z_\nu$ is the set of pairs $\Big((C,{\bf x},\ph),f\in F_{\nu,(C, {\bf x},
\ph)}\Big)$ obeying
\[
  \dbar\ph\ =\ \tau_\nu(f)\,.
\]
By construction, $\tau_\nu$ has support away from $U_0$. So $\ph$ is
holomorphic on $U_0(s)$ and the above equation can
actually be viewed on a local uniformizer as map
\[
  \dbar-\hat\tau: \calb\times \hat F_0\ \longrightarrow\ 
  \prod_{j\le m'} L^p(\Delta,\cz^n)\,,
\]
$\calb$ a complex Banach orbifold of the form given in Section~2.3. This map is
holomorphic with differential relative $S$ an epimorphism with finite
dimensional kernel. An application of the holomorphic implicit function theorem
with parameter space $S$ shows that $\tilde Z_\nu$ is locally uniformized by
a complex manifold. If $(C,{\bf x})$ is not stable we also have to take the
quotient by the germ of the action of $\Aut^0 (C,{\bf x})$. In \cite[\S5.3]{si1}
this has been achieved by imposing an averaged version (involving integrals
over bump functions in $M$) of the rigidification procedure given in
Proposition~\ref {CM_analyt_space}. The reason was that transversality
is not on open condition in $L_1^p$-spaces over two-dimensional domains.
On the subspace $\calb \subset S\times \check L_1^p(C;\ph^*T_M)$ we can
however use rigidification by incidence conditions with local transversal
divisors as in Proposition~\ref {CM_analyt_space}, as long as the
intersection is on the holomorphic part of $\ph$. The proof that this is
in fact a quotient is a simple application of the implicit function theorem,
parallel to the discussion in \cite[\S5.3]{si1}. Obviously, this slice by
incidence with transversal divisors is holomorphic, hence a complex
manifold. And the unrigidified $Z(q_\nu^*s_\dbar+\tau_\nu)$ is just a
product of the slice and an open set in $\cz^l$, $l=\dim \Aut(C,{\bf x})$.
\qed


\subsection{Reduction to local, holomorphic situation}
From the Kuranishi structure $\tau: F\rightarrow \cale$
(Proposition~\ref{tau}) we obtain a cone $C^F\subset F$ supporting a class
$[C^F]$ of dimension $d+\rk F$ as cone and class associated to the tautological
section of $q^*F|_{\tilde Z}$, $d=d(M,R,g,k)$ the expected dimension of
$\calc_{R,g,k} (M)$, $q:F \rightarrow \calc(M;p)$ the bundle projection, $\tilde Z
= Z(q^*s+\tau)$. Recall that the restriction of $F$ to $Z$ decomposes into a
direct sum of copies of $H$. We define
\[
  \mu:\ H\longhookrightarrow F|_Z
\]
to be the diagonal embedding. Let $\Theta _{F/H}$ be the Thom class of the
orbibundle $F/ \mu(H)$. Pulling back to $F$ yields a class $\Theta\in H_H^{\rk
F-\rk H}(F|_Z)$. Then $[C^F]\cap \Theta \in H_{d+\rk H} (C^F\cap H)$
is the intersection of $[C^H]$ with $H$. The comparison theorem
will readily follow from
\begin{prop}\label{main_prop}
  $C^H\ =\ C^F\cap H$, $[C^H]\ =\ [C^F]\cap\Theta$ in $H_{d+\rk H}(C^F\cap H)$.
\end{prop}
The proof of this proposition occupies the rest of this chapter. Notice that
since $C^H$ is a complex space of pure dimension $d+\rk H$, the homological
statement concerns classes of top dimension and can thus be checked locally
together with the set-theoretic part of the claim.
\vspace{2ex}

By Lemma~\ref{unif_cones} and the definition of $C^H$ it suffices to check the
claims of the proposition on the level of local uniformizers. Since
the local uniformizers in the complex analytic
(Proposition~\ref{CM_analyt_space}) and in the complex Banach manifold
treatment (Proposition~\ref{tau}) can be obtained by incidence conditions
with the same set of divisors $H_1,\ldots, H_l$, and all objects are just trivial
products of the restriction to the slice with an $l$-dimensional smooth
space, it even suffices to work on unrigidified charts. For $(C,{\bf x},\ph)
\in \calc_{R,g,k}(M)$ let $\tau_\nu$ be as in Proposition~\ref{tau}.
As usual we mark liftings to the local uniformizing systems
under study by a hat and with subscript $0$ the restriction to the
center of a local uniformizing system.
\begin{lemma}
  Possibly after shrinking $U$ to a smaller neighbourhood of $(C,{\bf x}, \ph)$
  there exists a topological decomposition
  \[
    \hat F\ =\ \bar F\oplus \hat F_\nu
  \]
  with $\tau(\bar F)\subset\img\sigma$ and such that
  \[
    \hat C^F\ =\ \bar F\oplus \hat C^h,\quad
    [\hat C^F]\ =\ [\bar F]\oplus [\hat C^h]\,,
  \]
  where $\hat C^h$, $[\hat C^h]$ are the cone and cone class obtained from
  $\hat \tau_\nu$.
\end{lemma}
\pf
By invoking Proposition~\ref{cone_decomp} we just have to define
$\bar F$. By the Fredholm property of $\sigma=\dbar$ and since $\hat\tau$
spans the cokernel of $\hat\sigma$ the family (over $S$) of linear maps
\[
  S\times\Big(F_0\oplus{\check L}_1^p(C;\ph^*T_M)\Big)\ \longrightarrow\ 
  {\check L}^p (C;\ph^*T_M)\,,\quad
  (s,f,v)\longmapsto \hat\tau_{s,v}(f)-\dbar v
\]
consists of split epimorphisms. An application of the implicit function
theorem thus shows that $\hat T_\tau:=\hat F\oplus_{\hat\cale} \hat\calt$
is a topological vector bundle on $\hat U$ of rank $\rk F +d$.
It fits into a quasi-isomorphism
\[\begin{array}{ccc}
  \hat T_\tau&\longrightarrow&\hat\calt\\[2pt]
  \diagl{\hat\rho}&&\diagr{\hat\sigma}\\[6pt]
  \hat F&\stackrel{\hat\tau}{\longrightarrow}&\hat\cale
\end{array}\]
Therefore $\hat\rho_0^{-1}(\hat F_\nu)$ is a linear subspace of
$\hat T_{\tau,0}$ of dimension $\rk F_\nu+d$. Let $P\subset \hat T_\tau$
be a subbundle restricting to a complementary subspace to
$\hat\rho_0^{-1}(\hat F_\nu)$ in $\hat T_{\tau,0}$. Then, possibly
after going over to a smaller local uniformizing
system, we may set $\bar F:=\hat \rho(P)$.
\qed
\begin{lemma}
  With the identification $H=F_\nu$ it holds
  \[
    \hat\mu^{-1}(C^F)\ =\ \hat C^h,\quad
    [C^F]\cap\Theta\ =\ [\hat C^h]\,.
  \]
  \vspace{-4ex}
  
\end{lemma}
\pf
Consider the family of morphisms
\[
  \hat\mu_t\ =\ (t\cdot\id,\ldots,t\cdot \id, \id, t\cdot \id,\ldots,t\cdot \id):
  H\longrightarrow F=\oplus_\lambda F_\lambda
\]
with ``$\id$'' at the $\nu$-th entry. By the previous lemma the claim holds with
$\hat\mu_0$ replacing $\hat\mu$. Since $\hat\mu_t$ is a proper homotopy
between $\hat\mu_0$ and $\hat\mu= \hat\mu_1$ we just have to show
\[
  \mu_t^{-1}(\hat C^F)\ =\ \hat C^h
\]
for any $t$. To verify this on the fiber over some
$\hat z\in \hat U$ let $R^F\subset \hat F_{\nu, \hat z}$
map isomorphically to $R:= \coker \hat\sigma_{\hat z}$. Another
application of Proposition~\ref{cone_decomp} with a larger complementary
space $\bar F$ shows
\[
  \hat C^F\cap\hat{F_{\hat z}}\ =\ q^{-1}(C^R)\,,
\]
for some cone $C^R\subset R$ where $q: R^F\rightarrow R$
is the quotient map. Letting $q^H: \hat H_{\hat z}
\rightarrow R^H$ be the cokernel of $\hat G_{\hat z} \rightarrow \hat H_{\hat
z}$ and $\bar\mu_t: R^H \rightarrow R$ be the map induced by $\mu_t$, we
obtain
\[
  \mu_t^{-1}(\hat C^F)\ =\ (q^H)^{-1}\bar\mu_t^{-1}(\hat C^R)\,.
\]
But $\bar\mu_t =\lambda_t\cdot \bar\mu_0$ for some $\lambda_t\in \rz_{>0}$,
for the maps from $R^H$ to $\coker(\dbar: \hat\calt_{\hat z} \rightarrow
\hat\cale_{\hat z})$ induced by any of the $\tau_\nu$ all coincide. Hence
\[
  \mu_t^{-1}(\hat C^F)\ =\ (q^H)^{-1}\bar\mu_0^{-1}(\hat C^R)\ 
  =\ \hat\mu_0^{-1}(\hat C^F)\ =\ \hat C^h\,.
\]
\vspace{-7ex}

\qed
\vspace{2ex}

To prove Proposition~\ref{main_prop} it remains to compare $\hat C^h$ and
$\hat{C^H}$, which will be the concern of the next section.


\subsection{Comparison of holomorphic normal spaces}
We consider the following situation: Let $q:\calc \rightarrow S$ be a prestable
curve over a smooth parameter space $S$ (this will be applied to an
analytically semiuniversal deformation of $(C,{\bf x})$), $Z=\Hom_S(\calc ,M)$
with universal curve and universal morphism
\[
  \pi:\Gamma\longrightarrow Z,\quad \Phi:\Gamma\rightarrow M\,.
\]
We embed $Z$ into a complex manifold $\tilde Z$ as in Proposition~\ref{tau}
(where the present $\tilde Z$ is denoted $\tilde Z_\nu$). Explicitely, possibly
after shrinking $S$, we work on the complex Banach orbifold $\calb$ of
$L_1^p$-maps from fibers $C_s$ of $q$ to $M$, holomorphic away from a union
of small disks $\bigcup \Delta^j$, as constructed in Section~2.3. According to
Theorem~\ref{hol_str}, $Z$ is the fiber over $0$ of the holomorphic Fredholm
map
\[
  \Theta:\ \calb\longrightarrow\prod_{j\le m'} L^p(\Delta^{1/2};\cz^n)=:\cale',
  \quad \psi\longmapsto \Big(\dbar\psi|_{\Delta^j}\Big)\,,
\]
in appropriate holomorphic coordinates on $M$ and $U_i$. $\tilde Z$ on the
other hand is obtained from a holomorphic morphism $\tau: \calb\times\cz^r
\rightarrow \cale'$, spanning the cokernel of $G$ at any $z\in Z$, as fiber over
$0$ of
\[
  \tilde \Theta:\ \calb\times\cz^r\longrightarrow\cale',\quad
  (\psi,a)\longmapsto \tau_\psi(a)+G(\psi)\,.
\]
Let $\tilde \pi: \tilde \Gamma \rightarrow \tilde Z$ be the universal curve over
$\tilde Z$ (this is holomorphic) and $\tilde\Phi: \tilde\Gamma\rightarrow M$ the
(usually non-holomorphic) evaluation map. The maps defined so far fit into the
following diagram
\[\begin{array}{ccccc}
  \Gamma&\stackrel{j}{\llonghookrightarrow}&\tilde\Gamma &
  \stackrel{\tilde\Phi}{\llongrightarrow}&M\\[3pt]
  \diagl{\pi}&&\diagr{\tilde\pi}\\[3pt]
  Z&\stackrel{i}{\llonghookrightarrow}&\tilde Z
\end{array}\]
$Z$ can also be viewed as zero locus of the tautological section $s_\can$ of the
(trivial) bundle $F= \tilde Z\times\cz^r \downarrow \tilde Z
\subset\calb\times \cz^r$ (in the notation of Proposition~\ref{tau}, $F$
corresponds to $q_\nu^*F_\nu$, $q_\nu: \tilde Z_\nu\rightarrow \calb$ the
projection); or, using the given trivialization, as fiber over $0$ of the
projection
\[
  \prj_2: \tilde Z\subset \calb\times\cz^r\ \longrightarrow\ \cz^r\,.
\]
To such a description belongs a global normal space for $Z$ as follows: let
$\calf = \calo(F^\vee)$. Evaluation at $s_\can$ yields an epimorphism
\[
  \calf\ \longrightarrow\ \cali
\]
to the ideal sheaf of $Z$ in $\tilde Z$. We define $\psi^{-1}$ to be the
composition with the map $\cali \rightarrow \cali/\cali^2$ to the conormal
sheaf. Put for $\psi^0$ the identity morphism on $\Omega_{\tilde Z/S}|_Z$ and
$d:\cali/ \cali^2 \rightarrow \Omega_{\tilde Z/S}|_Z$ the differential. Then
\[
  \psi^\bullet:\ [\calf\stackrel{d\circ\psi^{-1}}{\rightarrow}
  \Omega_{\tilde Z|S}|_Z]\ \longrightarrow\ [\cali/\cali^2\stackrel{d}{
  \rightarrow}\Omega_{\tilde Z/S}/_Z]=\tau_{\ge1}\call^\bullet_{Z|S}
\]
is a global normal space for $Z$ relative $S$. The corresponding cone
$C(\psi^\bullet)\subset F$ is by definition the image of the normal cone
$C_{Z|\tilde Z}$ of $Z$ in $\tilde Z$ under
the embedding $N_{Z|\tilde Z} \hookrightarrow F$. By
Proposition~\ref{analyt_cone}
\[
  C(\Psi^\bullet)\ =\ C(s_\can),\quad
  [C(\psi^\bullet)]\ =\ [C(s_\can)]\,,
\]
set-theoretically and in homology respectively. Note that in the notation of the
previous section $C(s_\can)$ and $[C(s_\can)]$ are (non-rigidified versions of)
$C^h$ and $[C^h]$.

On the other hand we have the global normal space
\[
  \ph^\bullet:\ [\calh\rightarrow\calg]\ \longrightarrow\ \call^\bullet_{Z|S}
\]
constructed from $\pi$, $\Phi$ (Section~1.2). By construction $\calh=\calf$ and
$C(\ph^\bullet)$ is the non-rigidified version of $C^H$. To prove
Proposition~\ref{main_prop} (non-rigidified) we have to show $C(\ph^\bullet)
=C(\psi^\bullet)$. Since global normal spaces depend only on the
(ray of) map induced in cohomology \cite[Thm.3.3]{si2} this will follow from
\begin{lemma}
  Locally, there exists an invertible holomorphic function $\chi$
  and an isomorphism
  \[
    \lambda^\bullet:\ [\calh\rightarrow\calg]\ \longrightarrow\ 
    [\calh\rightarrow\Omega_{\tilde Z/S}|_Z]
  \]
  with $\lambda^{-1}=\id$ and $H^i(\psi^\bullet\circ\lambda^\bullet) =
  \chi\cdot H^i(\ph^\bullet)$, $i=-1,0$.
\end{lemma}
The proof will occupy the rest of this section. The maps in the lemma fit into
the following diagram, which is claimed to commute in cohomology
up to multiplication by $\chi$.
\[\begin{array}{ccc}
[\calh\rightarrow\Omega_{\tilde Z/S}|_Z]\\[-3pt]
&{\scriptstyle\ \psi^\bullet}\\[-2pt]
&\arrowltrb\\[-8pt]
{\scriptstyle \lambda^\bullet\ }
\begin{picture}(0,0)\unitlength 1pt\put(0,-15){\vector(0,1){30}}\end{picture}
\phantom{\scriptstyle \lambda^\bullet\ }
&&\quad [\cali/\cali^2\rightarrow\Omega_{\tilde Z/S}|_Z]\\[-12pt]
&\arrowlbrt\\[3pt]
&{\scriptstyle\ \ph^\bullet}\\[-6pt]
[\calh\rightarrow\calg]
\end{array}\]

\noindent\nopagebreak
{\bf\boldmath Definition of $\lambda^\bullet$}\\[2ex]
To define $\lambda^0$ we observe that since $\tilde Z$
is solution to the equation $\dbar\ph=\tau_\nu(C,{\bf x},\ph)(h)$ the
tangent bundle $T_{\tilde Z/S}$ is canonically isomorphic to
$H\oplus_{\breve\cale}\breve\calt$. By Lemma~\ref{GH_quasiiso},
this fibered product is in turn canonically isomorphic to $G$.
We define $\lambda^0$ as the dual of the composition
\[
  T_{\tilde Z/S}|_Z\ \simeq\ H\oplus_{\breve\cale}\breve\calt\ \simeq\ G\,.
\]
It is clear from the construction that $\lambda^\bullet: [\calh\rightarrow
\calg] \rightarrow [\calh\rightarrow \Omega_{\tilde Z/S}|_Z]$ is a {\em
commutative} square.
\vspace{1cm}

\noindent\nopagebreak
{\bf\boldmath Replacing $[\calh\rightarrow\calg]$ by a
\v Cech complex}\\[2ex]
To begin with we simplify the problem by
dropping a common part from the definition of $\lambda^\bullet$
and $\ph^\bullet$ as follows. First note that $[\calh\rightarrow \calg]$
represents $[R\pi_*\Phi^*\calt_M]^\vee$ in the derived sense.
Derived objects such as $[R\pi_*\Phi^*\calt_M]^\vee$
are unique up to unique isomorphism in the derived category and the
existence of such an isomorphism is what we mean by ``represent''.
But the last steps in the construction of $\ph^\bullet$ consisted
of the composition of this isomorphism with duality
\[
  [\calh\rightarrow \calg]\ \simeq\ [R\pi_*\Phi^*\calt_M]^\vee
  \ \simeq\ R\pi_*(\Phi^*\Omega_M\otimes\omega)\,,
\]
so we may as well drop this composition and work with $R\pi_*
(\Phi^*\Omega_M\otimes\omega)$ directly. The latter in turn can be
represented by the \v Cech complex\
\[
  [\pi^{(0)}_*(\Phi^*\Omega_M\otimes\omega)
  \rightarrow\pi^{(1)}_*(\Phi^*\Omega_M\otimes\omega)]\,,
\]
and this gives an explicit identification of the cohomology
of $\calh \rightarrow\calg$ with $R^i\pi_*(\Phi^*\Omega_M \otimes
\omega)$, $i=0,1$.
\vspace{1ex}

\noindent
Similarly we may factor $\lambda^\bullet$ over the morphism of
complexes of holomorphic Banach bundles representing relative
duality (Proposition~\ref{duality})
\[\begin{array}{cccrcl}
  \pi_*^{(0)}(\Phi^*\Omega_M\otimes\omega)
  &\hspace{-1ex}\longrightarrow\hspace{-1ex}&
  (\pi_*^{(1)}\Phi^*\calt_M)^\vee
  &\quad(\alpha_i)&
  \hspace{-1ex}\longmapsto\hspace{-1ex}
  &\Big( (v_{ij})\mapsto \displaystyle \sum_{i,j}
  \int_{\Gamma/Z}\textstyle (\frac{\alpha_i+\alpha_j}{2})(v_{ij})
  \wedge\dbar\rho_i\Big)\\
  \longdownarrow&&\longdownarrow\\[3pt]
  \pi_*^{(1)}(\Phi^*\Omega_M\otimes\omega)
  &\hspace{-1ex}\longrightarrow\hspace{-1ex}&
  (\pi_*^{(0)}\Phi^*\calt_M)^\vee&
  (\alpha_{ij})&\hspace{-1ex}\longmapsto\hspace{-1ex}
  &\Big( (v_k)\mapsto \displaystyle \sum_{i,j}
  \int_{\Gamma/Z}\textstyle \alpha_{ij}(\frac{v_i+v_j}{2})
  \wedge\dbar\rho_i\Big)
\end{array}\]
In fact, this duality morphism is the composition of the (topological)
dual of 
\[
  [\beta^{(0)},\beta^{(1)}]:\ [\pi_*^{(0)} \Phi^* \calt_M
  \rightarrow \pi_*^{(1)}\Phi^* \calt_M]\
  \longrightarrow\ [\breve \calt\rightarrow \breve\cale]
\]
from Lemma~\ref{GH_quasiiso} (restricted to $Z$) and
\[\begin{array}{cccrcl}
  \pi_*^{(0)}(\Phi^*\Omega_M\otimes\omega)& \llongrightarrow&
  (\breve\cale)^\vee&\quad (\alpha_i)&\longmapsto&\Big( \gamma \mapsto
  \sum_i\int_{\Gamma/Z} \rho_i \alpha_i\wedge\gamma\Big)\\[3pt]
  \longdownarrow&&\longdownarrow\\[3pt]
  \pi_*^{(1)}(\Phi^*\Omega_M\otimes\omega)& \llongrightarrow&
  (\breve\calt)^\vee&\quad (\alpha_{ij})&\longmapsto&\Big( v \mapsto
  \sum_{i,j}\int_{\Gamma/Z} \alpha_{ij}(v)\wedge \dbar\rho_i\Big)\,.
\end{array}\]
To verify this one needs a little computation. The composition of the
two upper horizontal arrows applied to a local holomorphic section
$(\alpha_i)$ of $\pi_*^{(0)}(\Phi^*\Omega_M\otimes\omega)$ and
evaluated at a section $(v_{jk})$ of $\pi^{(1)}_*\Phi^*\calt_M$ leads to
the fiber integral
\[
  \frac{1}{2}\sum_{i,j,k} \int_{\Gamma/Z} \rho_i\alpha_i(v_{jk}\dbar \rho_j)\,.
\]
This indeed agrees with the upper horizontal arrow of the duality morphism
by noting that a partial integration computation shows
\[
  \int_{\Gamma/Z} \alpha_i(v_{ij}\dbar \rho_i)\ =\
  2\int\rho_i \alpha_i(v_{ij}\dbar \rho_i\,.
\]
Similarly for the lower horizontal arrows.

We may therefore draw a commutative diagram
\begin{eqnarray}\label{diag0}\begin{array}{ccccccccc}
  H^\vee&\hspace{-1ex}\longleftarrow\hspace{-1ex}&
  (\pi_*^{(1)}\Phi^*\calt_M)^\vee&
  \hspace{-1ex}\longleftarrow\hspace{-1ex}&
  \pi_*^{(0)}(\Phi^*\Omega_M\otimes\omega)&
  \hspace{-1ex}\longrightarrow\hspace{-1ex}&
  (\breve\cale)^\vee& \hspace{-1ex}\longrightarrow
  \hspace{-1ex}&H^\vee\\[3pt]
  \longdownarrow&&\longdownarrow&&\longdownarrow&&
  \longdownarrow&&\longdownarrow\\[3pt]
  G^\vee\hspace{-1ex}&\hspace{-1ex}\longleftarrow\hspace{-1ex}&
  (\pi_*^{(0)}\Phi^*\calt_M)^\vee&
  \hspace{-1ex}\longleftarrow\hspace{-1ex}&
  \mbox{\boldmath $\pi_*^{(1)}(\Phi^*\Omega_M\otimes\omega)$}&
  \hspace{-2ex}\mbox{\boldmath $\longrightarrow$}\hspace{-2ex}&
  \mbox{\boldmath $(\breve\calt)^\vee$}&
  \hspace{-2ex}\mbox{\boldmath $\longrightarrow$}\hspace{-1ex}&
  \mbox{\boldmath $T_{\tilde Z/S}^\vee|_Z$}\\[3pt]
  &&&&{\thicklines\longdownarrow}
  &&&&{\thicklines\longdownarrow}\\[3pt]
  &&&&\hspace{-2ex}\mbox{\boldmath $R^1\pi_*
  (\Phi^*\Omega_M\otimes\omega)$}\hspace{-2ex}&&
  \stackrel{H^0(\lambda^\bullet)}{\begin{picture}(0,6)(2.5,-4)
  \unitlength 1pt\put(-30,0){\vector(1,0){65}}\end{picture}}
  &&\mbox{\boldmath $\Omega_{Z/S}$}
\end{array}\hspace{3.5ex}\end{eqnarray}
The first three squares are compatible with the dual of $[G\rightarrow H]
\rightarrow [\breve\calt \rightarrow \breve \cale]$. It thus suffices to
prove the claim with $\lambda^\bullet$ replaced by the composition
of the right-hand squares in Diagram~\ref{diag0}.
\vspace{1cm}

\noindent\nopagebreak
{\bf\boldmath Computation of $H^0(\psi^\bullet\circ\lambda^\bullet)
= H^0(\lambda^\bullet)$}\\[2ex]
To describe $H^0(\lambda^\bullet)$ we now follow the bold printed part of the
diagram. Let $( \alpha_{ij}) \in \pi_*^{(1)}(
\Phi^* \Omega_M\otimes\omega)$ be a holomorphic family of \v Cech cochains
representing a local holomorphic section $\alpha$ of $R^1\pi_* 
(\Phi^*\Omega_M \otimes\omega)$. The associated local holomorphic section of
$(\breve\calt)^\vee$ is
\begin{eqnarray}\label{eqn4}
  {\check L}_1^p (C;\ph^*T_M)\cap \calo(\ph^*T_M|_{U_0})
  =\breve \calt_{(C,{\bf x},\ph)}\ni v\ \longmapsto\ \sum_{i,j}
  \int_{\Gamma/Z} \alpha_{ij}(v)\wedge \dbar\rho_i\,.
\end{eqnarray}
The restriction of this section to $Z\subset \calb$ is $H^0
(\lambda^\bullet)(\alpha)$.
\vspace{1cm}

\noindent
\nopagebreak
{\bf\boldmath Comparison with $H^0(\ph^\bullet)$}\\[2ex]
On the zeroth cohomology $\ph^\bullet$ is described by the
following diagram
\[\begin{array}{ccccc}
  \pi_*^{(1)}(\Phi^*\Omega_M \otimes\omega)&\llongrightarrow&
  \pi_*^{(1)}(\Omega_{\Gamma/\calc}\otimes\omega)&
  \stackrel{\simeq}{\llongrightarrow}&
  \pi_*^{(1)}(\pi^*\Omega_{Z/S}\otimes\omega)\\[3pt]
  \longdownarrow&&&&\longdownarrow\\[3pt]
  R^1\pi_*(\Phi^*\Omega_M \otimes\omega)&&&&\Omega_{Z/S}
\end{array}\]
The first horizontal arrow is by pull-back with $\Phi^*: \Phi^* \Omega_M
\rightarrow \Omega_\Gamma$ composed with the quotient $\Omega_\Gamma
\rightarrow \Omega_{\Gamma/\calc}$; the second horizontal arrow is by the
isomorphism $\Omega_{\Gamma/\calc} \simeq \pi^*\Omega_{Z/S}$; and the
right-hand vertical arrow is by projection formula composed with
the trace morphism $R^1\pi_* \omega \simeq \calo_Z$. Explicitely,
the map is (up to multiplication by an invertible holomorphic function
$\chi$ that we suppress here and in the sequel)
\begin{eqnarray}\label{eqn5}
  (\alpha_{ij}) &\longmapsto& \sum_{i,j} \int (\Phi^*
  \alpha_{ij})\wedge \dbar\rho_i\,,
\end{eqnarray}
interpreted as section of $\Omega_{Z/S}$ (that is, $\Phi^*$ taken as
pull-back of differential forms). Again we used the explicit version of relative
duality given in Proposition~\ref{duality}.
\begin{lemma}
  $H^0(\ph^\bullet)\ =\ H^0(\lambda^\bullet)$.
\end{lemma}
\pf
We evaluate $H^0(\ph^\bullet)(\alpha)$ at a local holomorphic section $v$
of $\breve \calt$. Since $\supp\dbar\rho_i \subset \bigcup_{j\le m'} \Delta^j$
we may restrict attention to $\alpha$ with only non-zero component
$\alpha_{ij}$. With local holomorphic coordinates $t$ on $\Delta$, $s$ on $S$
and $w^\mu$ on $M$ we may write
\[
  \alpha_{ij}\ =\ \sum_\mu a_\mu(t,s,\ph)d w^\mu\otimes dt\,,
  \quad v\ =\ \sum_\mu v^\mu(t,s, \ph)\partial_{w^\mu}
\]
for the relevant local parts. We obtain
\begin{eqnarray*}
  H^0(\ph^\bullet)(\alpha)(v)&=&
  \Big(\int_{\Gamma/Z} \Phi^*\alpha_{ij}\wedge\dbar\rho_i\Big)(v)
  \ =\ \int_\Delta \sum_\mu a_\mu(t,s,\ph)v^\mu (t,s,\ph)\, dt\wedge\dbar\rho_i\\
  &=&\int_\Delta \alpha_{ij}(v)\dbar\rho_i\ =\ H^0(\lambda^\bullet)(\alpha)(v)\,.
\end{eqnarray*}
\nopagebreak
\vspace{-8ex}

\qed
\vspace{1cm}
\pagebreak

\nopagebreak
\noindent
{\bf\boldmath Computation of $H^{-1}(\psi^\bullet\circ \lambda^\bullet)=
H^{-1}(\psi^\bullet)$}\\[2ex]
For holomorphic sections of $H^\vee$ coming from $(\breve\cale)^\vee$
(cf.\ Diagram~\ref{diag0}), pairing with the tautological
section of $H$ is nothing but pairing with the section $\breve s_\dbar$ of
$\breve \cale$. Given a local holomorphic section $\alpha$ of $\pi_*(\Phi^*
\Omega_M\otimes \omega)$ let $(\tilde \alpha_i)$ be a local holomorphic
section of $\tilde \pi_*^{(0)}(\tilde \Phi^*\Omega_M \otimes\omega)$
extending $\alpha|_{U_i}$. Then as map to the ideal sheaf $\cali$ of $Z$ in
$\tilde Z$,
\[
  H^{-1}(\psi^\bullet)(\alpha)\ =\ \sum_i\int \rho_i\tilde\alpha_i (\tilde
  s_\dbar)\ =\ \Big(\calb\ni (s,\ph)\mapsto \sum_i\int_{\tilde\Gamma/ \tilde Z}
  \rho_i\tilde \alpha_i(\dbar\ph)\Big)\Big|_Z\,.
\]
The induced section of $\cali/\cali^2$ depends only on $\alpha$, not on the
choice of extension $(\tilde\alpha_i)$.
\vspace{1cm}

\noindent
\nopagebreak
{\bf\boldmath Computation of $H^{-1}(\ph^\bullet)$}\\[2ex]
This step is a little harder and the most interesting part of the proof.
It will show that the $\dbar$-operator naturally turns up by an integration
by parts. Recall that $\ph^\bullet$ was defined in Section~1.2
by $R\pi_* (\,.\,\otimes\omega)$ of
\begin{eqnarray}\label{eqn6}
  L\Phi^*\call_M^\bullet\ \stackrel{L\Phi^*}{\longrightarrow}
  \call_{\Gamma/\calc}^\bullet\simeq L\pi^*\call_{Z/S}^\bullet\,.
\end{eqnarray}
All this is compatible with truncation $\tau_{\ge -1}$. Let $\cali\subset
\calo_{\tilde Z}$ and $\calj\subset \calo_{\tilde\Gamma}$ be the ideal sheaves
defining $Z\subset \tilde Z$ and $\Gamma\subset \tilde \Gamma$ respectively.
The following parts of the truncation of the previous sequence of complexes
are immediate:
\begin{eqnarray*}&&
  L\Phi^*\call_M^\bullet= [\Phi^*\Omega_M]\quad
  \mbox{\rm (one term in degree zero)}\\&&
  \tau_{\ge-1}\call_\Gamma^\bullet=[\calj/\calj^2\rightarrow \Omega_{
  \tilde\Gamma}|_\Gamma]\ \rightarrow\ \tau_{\ge-1}
  \call_{\Gamma/\calc}^\bullet=[\calj/\calj^2\rightarrow
  \Omega_{\tilde\Gamma/\calc}|_\Gamma]\\&&
  \tau_{\ge-1} \call_{Z/S}^\bullet=[\cali/\cali^2\rightarrow
  \Omega_{\tilde Z /S}|_Z]\\&&
  \tau_{\ge-1} L\pi^*\call_{Z/S}^\bullet=[\pi^*(\cali/\cali^2)\rightarrow
  \pi^* \Omega_{\tilde Z/S}]=[\calj/\calj^2\rightarrow
\Omega_{\tilde\Gamma/\calc}|_\Gamma]\,,
\end{eqnarray*}
the last line by flatness of $\pi$. To work out $L\Phi^*$ we decompose $\Phi$
as a closed embedding into a smooth space followed by a projection
\[
  \Phi:\ \Gamma\ \stackrel{(\iota,\Phi)}{\longhookrightarrow}\ 
  \tilde\Gamma\times M\ \stackrel{p}{\longrightarrow}\ M\,.
\]
Let $\tilde \calj$ be the ideal sheaf of $\Gamma$ in $\tilde \Gamma\times M$.
Writing $\tau_{\ge-1} \call_\Gamma^\bullet= [\tilde\calj/ \tilde\calj^2
\rightarrow \Omega_{\tilde\Gamma\times M}|_\Gamma]$,
\[
  \tau_{\ge-1} L\Phi^*:\ [0\rightarrow \Phi^*\Omega_M]\ \longrightarrow\ 
  [\tilde\calj/ \tilde\calj^2 \rightarrow
  \Omega_{\tilde\Gamma}|_\Gamma\oplus \Phi^*\Omega_M]
\]
is nothing but the inclusion on the degree zero term. To go further we need the
explicit quasi-isomorphism between the two representations of $\tau_{\ge-1}
\call_\Gamma^\bullet$ obtained from the embeddings into $\tilde\Gamma$
and $\tilde\Gamma\times M$ respectively: It is given by the
natural inclusion
\[
  [\tilde\calj/ \tilde\calj^2 \rightarrow
  \Omega_{\tilde\Gamma}|_\Gamma]\ \longrightarrow\ 
  [\tilde\calj/ \tilde\calj^2 \rightarrow
  \Omega_{\tilde\Gamma}|_\Gamma\oplus \Phi^*\Omega_M]\,.
\]
The truncation of sequence (\ref{eqn6}) can thus explicitely be written
\[\begin{array}{ccccccccl}
  0&\llongrightarrow&\tilde\calj/\tilde\calj^2&\llongleftarrow&
  \calj/\calj^2&=\!=\!=\!=&\calj/\calj^2&\simeq
  &\pi^*(\cali/\cali^2)\\[3pt]
  \longdownarrow&&\longdownarrow&&\longdownarrow&&
  \longdownarrow\\[6pt]
  \Phi^*\Omega_M&\llongrightarrow&\Omega_{\tilde\Gamma}|_\Gamma
  \oplus\Phi^*\Omega_M&\llongleftarrow&\Omega_{\tilde\Gamma}|_\Gamma&
  \llongrightarrow&\Omega_{\tilde\Gamma/\calc}|_\Gamma&\simeq
  &\pi^*(\Omega_{\tilde Z/S}|_Z)
\end{array}\]
To write down $R\pi_*(.\otimes\omega)$ of this diagram we abbreviate
\[
  {\rm Ev}:=\Phi^*\Omega_M\,,\quad
  J:=\calj/\calj^2\,,\quad
  \tilde J:=\tilde\calj/\tilde\calj^2\,,
\]
and, for a sheaf $\calf$ on $\Gamma$
\[
  \calf^{(\nu)}:=\pi_*^{(\nu)}(\calf\otimes\omega)\,.
\]
We also omit restrictions to $\Gamma$. For example,
\[
  (\Omega_{\tilde\Gamma/\calc}\oplus{\rm Ev})^{(1)}\ =\ 
  \pi_*^{(1)}\Big((\Omega_{\tilde\Gamma/\calc}|_\Gamma
  \oplus \Phi^*\Omega_M) \otimes\omega\Big)\,.
\]
A minor technical point arises when $\calf$ is not locally free. Then
$\calf^{(\nu)}$ is not a holomorphic Banach bundle over $Z$. Rather
there will be a holomorphic Banach bundle $E$ and a closed (ringed)
subspace $F^{(\nu)}\subset E$ given by finitely many holomorphic
functions that are linear in the fiber directions (i.e.\ $F^{(\nu)}$ is
a Banach version of a linear fiber space over $Z$ in complex analysis),
and $\calf^{(\nu)}$ is the sheaf of germs of holomorphic morphisms to the
trivial fiber space $Z\times\cz$. For our purposes the knowledge
of what a (holomorphic) section of $\calf^{(\nu)}$ is together with
the fact that kernel and cokernel of the \v Cech differential
$\calf^{(0)}\rightarrow \calf^{(1)}$ are the coherent sheaves
$\pi_*(\calf\otimes\omega)$ and $R^1\pi_*(\calf\otimes\omega)$
will suffice.

With these conventions the truncation of $R\pi_*(\,.\,\otimes\omega)$ applied
to sequence (\ref{eqn6}) is
\[\begin{array}{ccccccc}
  ({\rm Ev})^{(0)}&\hspace{-1ex}\longrightarrow\hspace{-1ex}&
  \Big((\Omega_{\tilde\Gamma/\calc}
  \oplus{\rm Ev})^{(0)}\oplus\tilde J^{(1)}\Big)\Big/ \tilde J^{(0)}&
  \hspace{-1ex} \longleftarrow\hspace{-1ex}
  &(\Omega_{\tilde\Gamma/\calc}^{(0)}\oplus J^{(1)})/ J^{(0)}
  &\hspace{-1ex}\longrightarrow\hspace{-1ex}&\cali/\cali^2\\[3pt]
  \longdownarrow&&\longdownarrow&&\longdownarrow
  &&\longdownarrow\\[6pt]
  ({\rm Ev})^{(1)}&\hspace{-1ex}\longrightarrow\hspace{-1ex}&
  (\Omega_{\tilde\Gamma/\calc}\oplus{\rm Ev})^{(1)}&
  \hspace{-1ex}\longleftarrow\hspace{-1ex}&
  \Omega_{\tilde\Gamma/\calc}^{(1)}&
  \hspace{-1ex}\longrightarrow\hspace{-1ex}&
  \Omega_{\tilde Z/S}|_Z
\end{array}\]
Here $\tilde J^{(0)}$, $J^{(0)}$ map to $(\Omega_{\tilde\Gamma/\calc}
\oplus{\rm Ev})^{(0)}$ and to $\Omega_{\tilde\Gamma/\calc}^{(0)}$
respectively by the K\"ahler differentials
\[
  d:\ \tilde\calj/\tilde\calj^2\ \longrightarrow\ \Omega_{\tilde\Gamma
  \times M}|_\Gamma \simeq \Omega_{\tilde\Gamma}|_\Gamma \oplus
  \Phi^*\Omega_M\,,\quad
  d:\calj/\calj^2\ \longrightarrow\ \Omega_{\tilde\Gamma}|_\Gamma\,,
\]
while the maps to $\tilde J^{(1)}$ and $J^{(1)}$ are by \v Cech differentials.
Similarly for the maps on \v Cech 1-cycles, but with one negative sign. The two
horizontal maps on the right are induced by the projection formula
composed with the trace morphism $R^1\pi_*\omega \simeq \calo_Z$
(Lemma~\ref{trace}).

Now let us chase some $\alpha\in \pi_*(\Phi^*\Omega_M \otimes\omega)$
through the upper horizontal sequence. Let $\alpha_i= \alpha|_{U_i}$. We claim
that $((0\oplus \alpha_i),0) \in ((\Omega_{\tilde\Gamma /\calc})^{(0)}
\oplus \tilde J^{(1)})$ lies in $\Omega_{\tilde\Gamma
/\calc}^{(0)} \oplus J^{(1)}$ modulo  $\tilde J^{(0)}$. From the exact sequence
\[
  0\ \longrightarrow\ \calj/\calj^2\ \longrightarrow\ 
  \tilde\calj/\tilde\calj^2\ \longrightarrow\ \Phi^*\Omega_M\
  \longrightarrow\ 0
\]
we obtain
\[
  0\ \longrightarrow\ J^{(\nu)}\ \stackrel{p_1^*}{\longrightarrow}\ 
  \tilde J^{(\nu)}\ \stackrel{d_M}{\longrightarrow}\ {\rm Ev}^{(\nu)}\
  \longrightarrow\ 0\,.
\]
There thus exists a local holomorphic
section $(\tilde g_i)$ of $(\tilde J)^{(0)}$ with
\[
  d_M \tilde g_i\ =\ -\alpha_i\,.
\]
But then, if $U_i\cap U_j\neq \emptyset$ (i.e.\ $i=0$ or $j=0$)
\[
  d_M(\tilde g_i-\tilde g_j)\ =\ \alpha_i-\alpha_j\ =\ 0\,,
\]
as section of ${\rm Ev}$. So the local section $(\tilde g_i-\tilde g_j)$ of
$\tilde J^{(1)}$ actually comes from a section $g_{ij}$ of $J^{(1)}$.
Therefore, $((0\oplus \alpha_i),0)$ lifts modulo $\tilde J^{(0)}$ to
the local section
\[
  (d_{\tilde \Gamma} \tilde g_i,g_{ij})
\]
of $\Omega_{\tilde\Gamma/\calc}^{(0)} \oplus J^{(1)}$. This maps to
$\cali/\cali^2$ by fiber integration
\begin{eqnarray}\label{eqn7}
  \frac{1}{2}\sum_{i,j}\int_{\Gamma/Z} g_{ij}\wedge\dbar\rho_i\,.
\end{eqnarray}
In view of the form of $H^{-1}(\psi^\bullet)$ discussed above
it remains to be shown:
\begin{lemma}\label{integral_descr}
  As local section of $\cali/\cali^2$ this integral equals
  \[
    \tilde Z\ni z\ \longmapsto\ \sum_i\int_{\tilde\Gamma_z}\rho_i
    \tilde\alpha_i(\dbar\ph_z)\,,
  \]
  where $(\tilde\alpha_i)$ extends $(\alpha_i)$ as section of $\tilde \pi_*
  (\tilde \Phi^*\Omega_M \otimes\omega)$, and we wrote $\ph_z:= \tilde
  \Phi|_{\tilde \Gamma/\tilde Z}$.
\end{lemma}
\pf
$(\id,\tilde\Phi): \tilde\Gamma \rightarrow \tilde\Gamma\times M$ is a (usually
non-holomorphic) extension of the closed embedding $\Gamma
\hookrightarrow \tilde\Gamma\times M$. The functions $(\id,\Phi)^*\tilde g_i$
will therefore not in general be holomorphic. But the fiber integrals
\begin{eqnarray}\label{eqn8}
  z\ \longmapsto\ \frac{1}{2}\sum_{i,j}\int_{\tilde\Gamma/\tilde Z}(\id,\Phi)^*
  (\tilde g_i-\tilde g_j) \wedge\dbar\rho_i\ =\ 
  \sum_i\int_{\tilde\Gamma/\tilde Z}(\id,\Phi)^* \tilde g_i \wedge\dbar\rho_i
\end{eqnarray}
will be local holomorphic functions on $\tilde Z$ as one easily sees in local
coordinates. Since $\tilde g_i-\tilde g_j$ induces the holomorphic
section $g_{ij}$ of $J^{(1)}$ the fact that $\dbar\rho_i=-\dbar \rho_j$
on $U_i\cap U_j$ shows that this holomorphic function induces the
section (\ref{eqn7}) of $\cali/\cali^2$.

On the other hand, partial integration applied to (\ref{eqn8}) results
in
\[
  -\sum_i\int_{\tilde\Gamma_z}\rho_i\dbar(\id,\Phi)^*\tilde g_i\ =\ 
  -\sum_i\int_{\tilde\Gamma_z}\rho_i d_M\tilde g_i(\dbar\ph_z)\,.
\]
Here we wrote $\ph_z=\Phi|_{\tilde\Gamma_z}$ and $d_M$ to denote the
composition
\[
  \calo_{\tilde\Gamma\times M}\ 
  \stackrel{d}{\longrightarrow}\ \Omega_{\tilde\Gamma\times M}\ 
  \simeq\ \Omega_{\tilde\Gamma}\boxplus\Omega_M
  \ \longrightarrow\\ p_2^*\Omega_M\,.
\]
Putting $\tilde\alpha_i:=-d_M \tilde g_i$ with $d_M\tilde g_i$ viewed as
holomorphic section of $\tilde\pi^{(0)}_*(\tilde\Phi^*\Omega_M
\otimes\omega)$ finishes the proof.
\qed


\section{Relative duality}
A (say algebraic) family of prestable curves $\pi:X\rightarrow S$ is Gorenstein:
It has an invertible relative dualizing sheaf $\omega_{X/S}$. For any coherent
sheaf $\calf$ on $X$ the theory of duality in derived categories developed in
\cite{hartshorne} then takes the following form. (We basically adopt the
terminology of loc.cit.\ except that we drop any underlining and we write
$\Ext$ for $R\,\Hom$.) It provides a {\em trace isomorphism}
\begin{eqnarray}\label{tracemap}
  \tr_\pi: R \pi_*\omega_{X/S}[1]\ \simeq\ R^1 \pi_*\omega_{X/S}\
  \longrightarrow\ \calo_S
\end{eqnarray}
and a {\em Yoneda morphism} (our notation)
\begin{eqnarray}\label{pairing}
  R \pi_* \Ext_X(\calf,\omega_{X/S})\ \longrightarrow\ 
  \Ext_S(R\pi_* \calf,R\pi_*\omega_{X/S})\,.
\end{eqnarray}
Composing we obtain the {\em duality morphism}
\begin{eqnarray}\label{duality_morphism}
  R \pi_* \Ext_X(\calf,\omega_{X/S})\ \longrightarrow\
  \Ext_S(R\pi_* \calf,\calo_S)= (R\pi_*\calf)^\vee\,,
\end{eqnarray}
and the content of the duality theorem is that this is an isomorphism.

The morphisms are in $D^+_{\rm coh}(S)$, the derived category of
complexes of $\calo_S$-modules bounded below
and with coherent cohomology. Sheaves are identified with
complexes concentrated in degree 0. For applicability of these
statements in the Gorenstein rather than the smooth case
see the remark at the beginning of \cite[VII.4]{hartshorne}.
A useful reference for all this is also \cite{lipman}. We will not
use duality theory for complex spaces, involving Fr\'echet sheaves
\cite{verdier}.

The purpose of this Chapter is to (a) make the transition from algebraic
sheaves to associated analytic sheaves (GAGA) and (b) to give an
explicit formulation of duality in terms of (analytic) \v Cech cochains and
fiber integrals, at least locally analytically over $S$.
\vspace{1ex}

Let us first comment on a subtlety that may be a source of confusion:
In our application it is crucial that all morphisms are unique up to unique
isomorphism (in the derived sense, algebraically or analytically).
Explicitely, this means that whenever we choose two representatives
of any of the objects (such as $R\pi_*\calf$) in terms of complexes of
$\calo_X$-modules there is a (sequence of) quasi-isomorphism(s)
between them that is unique up to homotopy. The same is true for
morphisms. In particular the maps induced in cohomology are
indeed unique up to unique isomorphism. In view of
\cite[Thm.3.3]{si2} this is enough to assure that the associated
(analytic or algebraic) cones are unique up to unique isomorphism and
hence are compatible with changes of local uniformizing systems.
\vspace{1ex}

Our plan is (1)~to express the Yoneda morphism in terms
of algebraic \v Cech-cochains (2)~to go over to analytic sheaves
and to admit refinements of the covering, and finally (3)~to give the
trace isomorphism (\ref{tracemap}) analytically. We do not indicate
in our notation if we are working algebraically or analytically.
For example, $X$ will denote either the scheme or its associated
analytic space, but the meaning will always be clear from the context.
\vspace{1ex}

Let $V_i$ be an affine open cover of $X$.
Then (the complex associated to) a
coherent sheaf $\calg$ on $X$ is quasi-isomorphic
to the complex of \v Cech-sheaves $[\calg^{(0)}\rightarrow
\calg^{(1)}]$. Recall that for any open $U\subset X$
the space of sections of $\calg^{(0)}$ over $U$ is
$\prod_i \calg(V_i\cap U)$, and similarly for $\calg^{(1)}$. The
corresponding \v Cech cochains relative $\pi$ are then just $\pi^{(i)}_*\calg
:= \pi_*\calg^{(i)}$.  Since the sheaves $\calg^{(i)}$ are $\pi_*$-acyclic
we can represent $R\pi_*\calg$ by $[\pi^{(0)}_*\calg\rightarrow
\pi^{(1)}_*\calg]$.

We claim that if $\calg$ is locally free then (possibly
after shrinking $S$) we can choose $V_i$ in such a way that $\pi^{(i)}_*
\calg$ are projective (i.e.\ locally free) $\calo_S$-modules.
To this end we factor $X\rightarrow S$ into a finite flat morphism
$\kappa: X\rightarrow S\times\pr^1$ and the projection $S\times \pr^1
\rightarrow \pr^1$. This is always possible after shrinking $S$.
We put $\bar V_0:= S\times(\pr^1\setminus\{\infty\})$,
$\bar V_1:= S\times(\pr^1\setminus\{0\})$ and $V_\nu=\kappa^{-1}(\bar
V_\nu)$. Denote the projections $\bar V_\nu\rightarrow S$ by $p_\nu$.
Now $\kappa_*\calg$ is locally free by flatness of $\kappa$. Writing
$S=\Spec A$, $\bar V_\nu=\Spec A[T]$ this means that the $A[T]$-module
$M$ associated to $\kappa_*\calg |_{\bar V_\nu}$ is projective. Equivalently,
it is a direct summand of a free $A[T]$-module. Viewed as an $A$-module,
$A[T]$ being a free $A$-module, $M$ is a direct summand of
a free $A$-module, hence projective. But the $A$-module $M$
is just the module associated to ${p_\nu}_* \kappa_*\calg$, so this
shows projectivity of $\pi^{(0)}_*\calg$. A similar argument
with $A[T]$ replaced by the ring of Laurent series $A[T]_{(T)}$
establishes the claim for $\pi^{(1)}_*\calg$.

We now restrict to locally free sheaves $\calf$ (slightly easier and
sufficient for our purposes).  By the digression we can then
represent $\Ext_S(R\pi_*\calf,
R\pi_*\omega_{X/S})$ by $\Hom^\bullet_S(\pi_*\calf^\bullet,
\pi_*\omega^\bullet)$. Written out this complex takes the form
\begin{eqnarray*}
  &&\Hom_S^\bullet\Big([\pi^{(0)}_*\calf\rightarrow \pi^{(1)}_*\calf],
  [\pi^{(0)}_*\omega_{X/S}\rightarrow \pi^{(1)}_*\omega_{X/S}]\Big)\\
  &=&\Big[\Hom_S(\pi^{(0)}_*\calf, \pi^{(0)}_*\omega_{X/S})
  \times\Hom_S(\pi^{(1)}_*\calf, \pi^{(1)}_*\omega_{X/S})\ \rightarrow\
  \Hom_S(\pi^{(0)}_*\calf, \pi^{(1)}_*\omega_{X/S})\Big]\,.
\end{eqnarray*}

\begin{lemma}\label{pairingmorphism}
  The Yoneda morphism $R\pi_*[\calf^\vee\otimes\omega_{X/S}]
  \rightarrow \Ext_S(R\pi_*\calf, R \pi_*\omega_{X/S})$ can be represented
  by the pair of morphisms $(\ph^0,\ph^1)$ with
  \begin{eqnarray*}
    \ph^0:\pi^{(0)}_*(\calf^\vee\otimes\omega_{X/S})&\longrightarrow&
    \Hom(\pi^{(0)}_*\calf, \pi^{(0)}_*\omega_{X/S})
    \times\Hom(\pi^{(1)}_*\calf, \pi^{(1)}_*\omega_{X/S})\\
    (\alpha_i)_i&\longmapsto&
    \Big((\beta_i)_i\mapsto (\alpha_i(\beta_i))_i,
    (\beta_{ij})_{ij}\mapsto \Big(\textstyle (\frac{\alpha_i+\alpha_j}{2})
    (\beta_{ij})\Big)_{ij}\Big)
  \end{eqnarray*}
  and
  \begin{eqnarray*}
    \ph^1:\pi^{(1)}_*(\calf^\vee\otimes\omega_{X/S})&\longrightarrow&
    \Hom(\pi^{(0)}_*\calf, \pi^{(1)}_*\omega_{X/S})\\
    (\alpha_{ij})_{ij}&\longmapsto&
    \Big((\beta_i)_i\mapsto \Big(\alpha_{ij}
    (\textstyle\frac{\beta_i+\beta_j}{2})\Big)_{ij}\Big)\,.
  \end{eqnarray*}
\end{lemma}
\pf
The Yoneda morphism is defined in \cite[II.5.5]{hartshorne}. The recipe
is to first represent $\calf$ and $\omega_{X/S}$ by a $\pi_*$-acyclic
resolution $\calf^\bullet$ and by an injective resolution $\omega^\bullet$
respectively. Then $\Hom_X^\bullet(\calf^\bullet, \omega^\bullet)$
represents $\Ext_X([\calf],[\omega])$ and consists of flasque sheaves.
Thus $R\pi_*\Ext_X([\calf],[\omega])$ can be represented by
$\pi_*\Hom_X^\bullet (\calf^\bullet, \omega^\bullet)$. Similarly, $\Ext_S
(R\pi_*[\calf],$ $R\pi_*[\omega])$ is represented by $\Ext_S
(\pi_*\calf^\bullet,\pi_*\omega^\bullet)$. With these representatives the
Yoneda morphism is simply the composition of natural morphisms
\[
  \pi_*\Hom_X^\bullet(\calf^\bullet, \omega^\bullet)\
  \longrightarrow\ \Hom_S^\bullet(\pi_*\calf^\bullet,\pi_*\omega^\bullet)
  \longrightarrow \Ext_S (\pi_*\calf^\bullet,\pi_*\omega^\bullet)\,.
\]
For $\calf^\bullet$ we may take the \v Cech complex $\calf^{(0)}\rightarrow
\calf^{(1)}$. And for our special choice of affine covering the
second arrow becomes an isomorphism and is thus understood.
We claim that we may take for $\omega^\bullet$ the \v Cech
resolution $\omega^{(0)}\rightarrow \omega^{(1)}$ as well, instead of an
injective one. In fact, $\Hom_X(\calf^{(i)},\omega^{(i)})$ consists of a direct
sum of coherent sheaves supported on affine sets and is thus
$\pi_*$-acyclic. By injectivity of $\omega^\bullet$ there exists a map
of complexes $[\omega^{(0)}\rightarrow \omega^{(1)}] \rightarrow
\omega^\bullet$. This map induces a commutative diagram of complexes
\[\begin{array}{ccc}
  \pi_*\Hom_X^\bullet([\calf^{(0)}\rightarrow \calf^{(1)}],
  [\omega^{(0)}\rightarrow \omega^{(1)}]) &\longrightarrow&
  \Hom_S^\bullet([\pi^{(0)}_*\calf\rightarrow \pi_*^{(1)}\calf],
  [\pi_*^{(0)}\omega \rightarrow\pi_*^{(1)}\omega])\\
  \longdownarrow&&\longdownarrow\\[6pt]
  \pi_*\Hom_X^\bullet([\calf^{(0)}\rightarrow \calf^{(1)}],
  [\omega^0\rightarrow\omega^1) &\longrightarrow&
  \Hom_S^\bullet([\pi^{(0)}_*\calf\rightarrow \pi_*^{(1)}\calf],
  [\pi_*\omega^0 \rightarrow\pi_*\omega^1])
\end{array}\]
in which {\em the vertical morphisms are quasi-isomorphisms}. End
proof of claim.

It remains to replace the resolution $\Hom_X^\bullet (\calf^\bullet,
\omega^\bullet)$ of $\Hom_X(\calf,\omega)$  by its \v Cech resolution
$\Hom_X(\calf,\omega)^{(0)}\rightarrow \Hom_X(\calf,\omega)^{(1)}$.
Explicitely, we have the following quasi-isomor\-phism
\[\begin{array}{ccccccc}
  0&\longrightarrow&\Hom(\calf,\omega)&\longrightarrow&
  \Hom^{(0)}(\calf,\omega)&
  \stackrel{\check d_{\rm Hom}} {\longrightarrow}&
  \Hom^{(1)}(\calf,\omega)\\
  &&\diagl{\rm Id}&&\diagr{F}&&\diagr{G}\\[6pt]
  0&\longrightarrow&\Hom(\calf,\omega)&\longrightarrow&
  \begin{array}{c}\Hom(\calf^{(0)},\omega^{(0)})\\[-3pt]
  \oplus\\ \Hom(\calf^{(1)},\omega^{(1)})\end{array}&
  \stackrel{\delta} {\longrightarrow}&
  \Hom(\calf^{(0)},\omega^{(1)})  
\end{array}\]
The maps are $\delta(a;b) = \check d_\omega\circ a - b\circ \check
d_\calf$,
\begin{eqnarray*}
  F:\ (\ph_i)_i&\longmapsto&\Big((f_i)_i\mapsto(\ph_i(f_i))_i\,;
  (f_{ij})_{ij}\mapsto \textstyle (\frac{\ph_i+\ph_j}{2})(f_{ij})\Big)\\
  G:\ (\ph_{ij})_{ij}&\longmapsto& \Big( (f_i)_i\mapsto
  \ph_{ij}(\textstyle\frac{f_i+f_j}{2})_{ij}\Big)\,,
\end{eqnarray*}
and $\check d_{\rm Hom}$, $\check d_\omega$, $\check d_\calf$ are
\v Cech differentials. Composing these maps with the
natural map to $\Hom_S^\bullet (\pi_*\calf^\bullet,
\pi_*\omega^\bullet)$ from above gives the stated result.
\qed
\vspace{1ex}

\noindent
For the next step (2) we view $X\rightarrow S$ as a morphism of complex
spaces. The analytic sheaves associated to the algebraic sheaves
$\pi_*^{(i)}(\calf^\vee\otimes\omega_{X/S})$ etc.\ are given by pull-back
under the morphism of ringed spaces from the complex space
$(S^{\rm an},\calo_{S^{\rm an}})$ underlying $S$ to the scheme $(S,\calo_S)$.
This amounts to going over to analytic topology and tensoring with
the structure sheaf of $S^{\rm an}$ over the pull-back of $\calo_S$. The
effect is that one considers analytic \v Cech cochains that are algebraic
along the fibers of $\pi$.

Since $\calf$ was supposed to be locally free the dualization of a sequence
analogous to (\ref{twist_sequence}) provides a short exact sequence of
algebraic sheaves on $X$
\begin{eqnarray}\label{free_res}
  0\ \longrightarrow\ \calf\ \longrightarrow\ \calg\ \longrightarrow\ \calh
  \ \longrightarrow\ 0
\end{eqnarray}
with $R^1\pi_*\calg=R^1\pi_*\calh=0$. Taking the associated $\pi_*$-acyclic
resolution by (algebraic) \v Cech sheaves $0\rightarrow \calf\rightarrow
\calf^{(0)} \rightarrow \calf^{(1)}\rightarrow 0$ etc.\ and pushing forward
by $\pi$ yields a diagram of a form (dual to) Diagram~\ref{big_diag}.
An argument similar to the one given there produces a
quasi-isomorphism (of quasi-coherent, algebraic sheaves)
\begin{eqnarray}\label{res_quis}
  [\pi_*\calg\rightarrow\pi_*\calh]\ \longrightarrow\ 
  [\pi_*^{(0)}\calf\rightarrow \pi_*^{(1)}\calf]\,.
\end{eqnarray}

For any point $s\in S$ we now choose a Stein refinement
$\{U_i\}$ of the affine covering $\{V_i\}$ of the form given in Section~2.3
(we assume that the point $0,\infty\in \pr^1$ ahev been chosen
suitably). Possibly by enlarging the index set for $\{V_i\}$ we may
assume that $\{U_i\}$ is indeed a shrinking of $\{V_i\}$ (in the
Hausdorff-topology). Now the analytic sheaves $\calf^{\rm an}$
etc. associated to $\calf$, $\calg$ and $\calh$ and the covering
$\{U_i\}$ give rise to a similar quasi-isomorphism (of
$\calo_{S^{\rm an}}$-modules)
\begin{eqnarray}\label{an_res_quis}
  [\pi_*\calg^{\rm an}\rightarrow\pi_*\calh^{\rm an}]\ \longrightarrow\ 
  [\pi_*^{(0)}\calf^{\rm an}\rightarrow \pi_*^{(1)}\calf^{\rm an}]\,.
\end{eqnarray}
We then obtain a restriction map from $\pi_*^{(i)}\calf\otimes
\calo_{S^{\rm an}}$ to the corresponding $\calo_{S^{\rm an}}$-module
of analytic \v Cech cochains associated to $\{U_i\}$. Since $\pi$
is a projective morphism, by Chows lemma
sections of an analytic sheaf $\call$ over $\pi^{-1}(U)$ are fiberwise
algebraic. In other words, $\pi_*\calg\otimes \calo_{S^{\rm an}}
\simeq \pi_*\calg^{\rm an}$ and similarly for $\calh$. Compatibility
of (\ref{an_res_quis}) with (\ref{res_quis}) tensored by $\calo_{S^{\rm an}}$
and the restriction map now shows the required GAGA statement
(just in this lemma we use for clarity the notation $\pi_*^{(i),{\rm an}}$
to denote relative bounded \v Cech cochains with respect to $\{U_i\}$ as
opposed to $\pi_*^{(i)}$ for fiberwise algebraic \v Cech cochains with respect
to $\{V_i\}$):
\begin{lemma}
  The restriction morphism
  \[
    [\pi_*^{(0)}\calf\otimes\calo_{S^{\rm an}}
    \rightarrow \pi_*^{(1)}\calf\otimes\calo_{S^{\rm an}}]\
    \longrightarrow
     [\pi_*^{(0),{\rm an}}\calf^{\rm an}\rightarrow
     \pi_*^{(1),{\rm an}}\calf^{\rm an}]
  \]
  is a quasi-isomorphism.
\qed
\end{lemma}
\vspace{1ex}

\noindent
As last ingredient (3) we give an analytic trace isomorphism. Let $\rho_i$ be
a partition of unity subordinate to $U_i$.
\begin{lemma}\label{trace}
  The map
  \[
    \Phi:\pi_*^{(1)}\omega_{X/S}\ \longrightarrow\ \calo_S\,,\quad
    (\alpha_{ij})_{ij}\ \longmapsto\ \sum_{i,j}\int_{X/S}\alpha_{ij}
    \wedge\dbar \rho_i
  \]
  induces an isomorphism $R^1\pi_*\omega_{X/S}\simeq\calo_S$.
\end{lemma}
\pf
The map $\Phi$ vanishes on the image of $\pi_*^{(0)}\omega_{X/S}$ under
the \v Cech differential and thus induces the claimed map on the
first cohomology. By duality on the fibers $X_s$ of $\pi$ it holds
$H^1(X_s, \omega_{X_s})= H^0(X_s,\calo_{X_s}) =\cz$ since $X_s$ (being
prestable) is reduced and connected. Naturality of relative dualizing
sheaves shows $\omega_{X_s} =\omega_{X/S}|_{X_s}$. The map $\pi$
being flat, Grauert's base change theorem thus implies that
$R^1\pi_*\omega_{X/S}$ is locally free of rank one with fibers
\[
  R^1\pi_*\omega_{X/S}/{\frak m}_s R^1\pi_*\omega_{X/S}\ 
  \simeq\ H^1(X_s,\omega_{X_s})\,.
\]
It thus suffices to find, for any prestable curve $C$, a 1-cocycle with values
in $\omega_C$ with non-vanishing value under the integral defining $\Phi$.
And indeed, any $U_{ij}\subset C$ being an annulus $\Delta\setminus
\bar\Delta_a$ (where we assume $U_i$ meet the inner boundary $\partial
\Delta_a$) we may put $\alpha_{ij}= z^{-1}dz$ with $z$ the linear
coordinate on $\Delta$. Then using integration by parts and
$\rho_i|_{\partial \Delta} \equiv 0$, $\rho_i|_{\partial\Delta_a}\equiv 1$
we obtain
\[
  \int_C\alpha_{ij}\wedge \dbar\rho_i\ =\
  \int_{\Delta\setminus\Delta_a} z^{-1} dz\wedge\dbar\rho_i\ =\ 
  \int_{\partial\Delta_a} z^{-1} dz\ =\ 2\pi i\ \neq\ 0\,.
\]
\vspace{-6ex}

\qed
\vspace{1ex}

Notice that an isomorphism $R^1\pi_*\omega_{X/S}\simeq\calo_S$
is unique up to multiplication by an invertible holomorphic function.
Thus the isomorphism given here differs from the analytification of
the algebraic trace isomorphism only by such multiplication.

To complete the explicit (analytic) description of algebraic relative
duality we just have to compose the trace isomorphism with the
Yoneda pairing.
\begin{prop}\label{duality}
  The algebraic duality morphism (\ref{duality_morphism}) is locally
  analytically given by
  \[\begin{array}{cccrcl}
  \pi_*^{(0)}(\calf^\vee\otimes\omega_{X/S})
  &\hspace{-1ex}\longrightarrow\hspace{-1ex}&
  (\pi_*^{(1)}\calf)^\vee
  &\quad(\alpha_i)&\hspace{-1ex}\longmapsto\hspace{-1ex}&
  \Big( (\beta_{ij})\mapsto \displaystyle \chi\cdot\sum_{i,j}
  \int_{\Gamma/Z}(\textstyle\frac{\alpha_i+\alpha_j}{2})(\beta_{ij})
  \wedge\dbar\rho_i\Big)\\[3pt]
  \longdownarrow&&\longdownarrow\\[3pt]
  \pi_*^{(1)}(\calf^\vee\otimes\omega_{X/S})
  &\hspace{-1ex}\longrightarrow\hspace{-1ex}&
  (\pi_*^{(0)}\calf)^\vee&
  (\alpha_{ij})&\hspace{-1ex}\longmapsto\hspace{-1ex}&
  \Big( (\beta_i)\mapsto \displaystyle \chi\cdot\sum_{i,j}
  \int_{\Gamma/Z}\alpha_{ij}(\textstyle\frac{\beta_i+\beta_j}{2})
  \wedge \dbar\rho_i\Big)
\end{array}\]
for some invertible holomorphic function $\chi$ on $S$.
\qed
\end{prop}


\noindent
{\small Department of Mathematics\\
Massachusetts Institute of Technology\\
Cambridge, MA 02139}\\[1ex]
{\small and}\\[1ex]
{\small Fakult\"at f\"ur Mathematik\\
Ruhr-Universit\"at Bochum, D-44780 Bochum\\
e-mail: Bernd.Siebert@ruhr-uni-bochum.de}
\end{document}